\documentclass[titlepage,draft,12pt]{article} 
\usepackage{amssymb,amsthm,amsmath} 

\usepackage[text={15.7 cm, 21.9 cm},centering,includefoot]{geometry}

\date{}

\newcommand{\ep}{\varepsilon}
\renewcommand{\qed}{{\penalty 10000\mbox{$\quad\Box$}}}
\newcommand{\re}{\mathbb{R}}

\newcommand{\n}{\mathbb{N}}

\newcommand{\rk}{\rho_{k}}
\newcommand{\tk}{\theta_{k}}
\newcommand{\ri}{\rho_{i}}
\newcommand{\ti}{\theta_{i}}

\newtheorem{thm}{Theorem}[section]

\newtheorem{rmk}[thm]{Remark}
\newtheorem{prop}[thm]{Proposition}
\newtheorem{defn}[thm]{Definition}

\newtheorem{lemma}[thm]{Lemma}

\title{Quantization of energy and weakly turbulent profiles of the solutions to some damped second order evolution equations}

\author{Marina Ghisi\vspace{1ex}\\ 
{\normalsize Universit\`a degli Studi di Pisa} \\
{\normalsize Dipartimento di Matematica}\\ 
{\normalsize PISA (Italy)}\\
{\normalsize e-mail: \texttt{marina.ghisi@unipi.it}}
\and
Massimo Gobbino\vspace{1ex}\\ 
{\normalsize Universit\`a degli Studi di Pisa} \\
{\normalsize Dipartimento di Matematica}\\ 
{\normalsize PISA (Italy)}\\  
{\normalsize e-mail: \texttt{massimo.gobbino@unipi.it}}
\and
Alain Haraux\vspace{1ex}\\ 
{\normalsize Universit\'{e} Pierre et Marie Curie} \\
{\normalsize Laboratoire Jacques-Louis Lions}\\ 
{\normalsize PARIS (France)}\\  
{\normalsize e-mail: \texttt{haraux@ann.jussieu.fr}}}

\begin{document}
\maketitle
\begin{abstract}

	We consider a second order equation with a linear ``elastic'' 
	part and a nonlinear damping term depending on a power of the 
	norm of the velocity. We investigate the asymptotic behavior of 
	solutions, after rescaling them suitably in order to take into 
	account the decay rate and bound their energy away from zero.
	
	We find a rather unexpected dichotomy phenomenon.  Solutions with finitely many Fourier
	components are asymptotic to solutions of the linearized equation
	without damping, and exhibit some sort of equipartition of the
	energy among the components.  Solutions with infinitely many
	Fourier components tend to zero weakly but not strongly.  We show
	also that the limit of the energy of solutions depends only on the
	number of their Fourier components.
	
	The proof of our results is inspired by the analysis of a
	simplified model which we devise through an averaging procedure,
	and whose solutions exhibit the same asymptotic properties as the
	solutions to the original equation.
	
\vspace{6ex}

\noindent{\bf Mathematics Subject Classification 2010 (MSC2010):} 
35B40, 35L70, 35B36.


\vspace{6ex}


\end{abstract}

 
\section{Introduction}

Let $H$ be a real Hilbert space, in which $|x|$ denotes the norm of an
element $x\in H$, and $\langle x,y\rangle$ denotes the scalar product
of two elements $x$ and $y$.  Let $A$ be a self-adjoint operator on
$H$ with dense domain $D(A)$.  We assume that $H$ admits a countable
orthonormal basis made by eigenvectors of $A$ corresponding to an
increasing sequence of positive eigenvalues $\lambda_{k}^{2}$.

We consider the second order evolution equation
\begin{equation}
	u''(t)+|u'(t)|^{2}u'(t)+Au(t)=0,
	\label{eqn:main}
\end{equation}
with initial conditions
\begin{equation}
	u(0)=u_{0}\in D(A^{1/2}),
	\quad\quad
	u'(0)=u_{1}\in H.
	\label{eqn:data}
\end{equation}

All nonzero solutions to (\ref{eqn:main}) decay to zero in the energy
space $D(A^{1/2})\times H$, with a decay rate proportional to
$t^{-1/2}$ (see Proposition~\ref{prop:basic}).  This suggests the
introduction and the investigation of the rescaled variable
$v(t):=\sqrt{t}\cdot u(t)$.\bigskip

The special structure of the damping term guarantees that for any linear subspace $F\subseteq D(A)$ 
such that  $A(F) \subseteq F$, the space $F\times F$ is positively invariant by  the flow generated by (\ref{eqn:main}).
In particular, equation (\ref{eqn:main}) possesses the so-called finite dimensional modes, namely solutions whose both components of the
initial state $(u_{0},u_{1})$ are finite combinations of the eigenvectors. Denoting by $u_{k}(t)$ and $v_{k}(t)$ the projections of $u(t)$ and $v(t)$ on 
the $k^{th}$ eigenspace, we shall call for simplicity the quantity 
$$ t\left(|u_{k}'(t)|^{2}+\lambda_{k}^{2}|u_{k}(t)|^{2}\right)$$ the ``energy of the $k^{th}$ Fourier component of $v(t)$" while 
$$ t\left(|u'(t)|^{2}+|A^{1/2}u(t)|^{2}\right)$$ will be called the ``energy of $v(t)$". For $t$ large, these quantities are easily seen to be equivalent
 to $ |v_{k}'(t)|^{2}+\lambda_{k}^{2}|v_{k}(t)|^{2}$ and $|v'(t)|^{2}+|A^{1/2}v(t)|^{2}$, respectively. Our main results, formally stated as Theorem~\ref{thm:main} and Theorem~\ref{thm:energ.}, can be
summed up as follows.
\begin{itemize}
	\item The limit of the energy of $v(t)$ depends only on the number
	of Fourier components of $v(t)$ that are different from~0.  In
	particular, the limit of the energy can take only countably many
	values.

	\item If $v(t)$ has only a finite number of Fourier components
	different from~0, then $v(t)$ is asymptotic in a strong sense to a suitable
	solution $v_{\infty}(t)$ to the nondissipative linear equation
	\begin{equation}
		v''(t)+Av(t)=0.
		\label{eqn:lin-hom}
	\end{equation}
	
	Moreover, there is \emph{equipartition of the energy} in the
	limit, in the sense that all nonzero Fourier components of
	$v_{\infty}(t)$ do have the same energy.

	\item If $v(t)$ has infinitely many components different from~0,
	then $v(t)$ tends to~0 weakly in the energy space, but not
	strongly.  Roughly speaking, the energy of $v(t)$ does not tend
	to~0, but in the limit there is again equipartition of the energy,
	now among infinitely many components, and this forces all
	components of $v(t)$ to vanish in the limit.
\end{itemize}

In other words, the Fourier components of rescaled solutions to
(\ref{eqn:main}) communicate to each other, and this can result in
some sort of energy transfer from lower to higher frequencies, longing
for a uniform distribution of the energy among components. In the case of an infinite 
number of non-trivial Fourier components, the weak convergence to $0$ implies 
non-compactness of the profile in the energy space. In particular, if $A$ has compact resolvent,
whenever the initial state $(u_{0},u_{1})$ belongs to $D(A)\times D(A^{1/2})$ and has an infinite number 
of elementary modes, the norm of $(v(t), v'(t))$ in $D(A)\times D(A^{1/2})$ is unbounded, a typical phenomenon 
usually called {\it weak turbulence}, cf.\ e.g.\ \cite{Bourg} and \cite{Open} for other examples. \bigskip

Our abstract theory applies for example to wave equations with 
nonlinear \emph{nonlocal} damping terms of the form
\begin{equation}
	u_{tt}(t,x)+
	\left(\int_{0}^{\ell}u_{t}^{2}(t,x)\,dx\right)u_{t}(t,x)
	-u_{xx}(t,x)=0
	\label{PDE:nonlocal}
\end{equation}
in a bounded interval $(0,\ell)$ of the real line with homogeneous Dirichlet boundary conditions. This is a toy 
model of the wave equation with \emph{local} nonlinear damping 
\begin{equation}
	u_{tt}(t,x)+u_{t}^{3}(t,x)-u_{xx}(t,x)=0
	\quad\quad
	t\geq 0,\ x\in(0,\ell),
	\label{PDE:local}
\end{equation}
which in turn is the prototype of all wave equations with nonlinear 
dissipation of order higher than one at the origin. This more general problem was the motivation that led us to consider
equations (\ref{PDE:nonlocal}) and (\ref{eqn:main}).   It is quite easy to prove
that all solutions to (\ref{PDE:local}) decay \emph{at least} as $t^{-1/2}$.  Actually the more general problem 
$$u_{tt}(t,x)+g(u_{t}(t,x))-\Delta u(t,x)=0$$ in any bounded domain with homogeneous Dirichlet boundary conditions and $g$ non-decreasing has been extensively studied under relevant assumptions on the behavior of $g$ near the origin and some conditions on the growth of $g$ at infinity, cf.\ e.g.\ \cite{Nakao, H-ondes, Dieud, HZ} in which reasonable energy estimates, of the same form as those in the ODE case are obtained. However, the
asymptotic behavior of solutions to the simple equation (\ref{PDE:local}) is still a
widely open problem since, unlike the ODE case, the optimality of this decay rate in unknown: there are neither examples of solutions to (\ref{PDE:local}) whose
decay rate is proportional to $t^{-1/2}$, nor examples of nonzero
solutions that decay faster.\bigskip

It is not clear whether our results shed some light on the local case
or not.  For sure, they confirm the complexity of the problem.  In the
case of (\ref{PDE:local}) there are no simple invariant subspaces, and
the interplay between components induced by the nonlinearity is more
involved.  Therefore, it is reasonable to guess that at most the
infinite dimensional behavior of (\ref{PDE:nonlocal}) extends to
(\ref{PDE:local}), and this behavior is characterized by lack of an
asymptotic profile and lack of strong convergence.\bigskip

As a matter of fact, the problem of optimal decay rates is strongly related to regularity
issues.  It can be easily shown that solutions to (\ref{PDE:local})
with initial data in the energy space remain in the same space for all
times, and their energy is bounded by the initial energy.  But what
about more regular solutions?  Can one bound higher order Sobolev
norms of solutions in terms of the corresponding norms of initial
data?  This is another open problem, whose answer would imply partial
results for decay rates, as explained in \cite{Below} and \cite{Open}, cf. also \cite{Van-Mart} for a partial optimality result in the case of boundary damping.  However, the
energy traveling toward higher frequencies might prevent the bounds
on higher order norms from being true, or at least from being easy to
prove.\bigskip

This paper is organized as follows.  In section~\ref{sec:statements}
we state our main results.  In section~\ref{sec:basic} we prove the
basic energy estimate from above and from below for solution to
(\ref{eqn:main}), we introduce Fourier components, and we interpret
(\ref{eqn:main}) as a system of infinitely many ordinary differential
equations.  In section~\ref{sec:heuristics} we consider a simplified
system, obtained from the original one by averaging some oscillating
terms.  Then we analyze the simplified system, we discover that it
is the gradient flow of a quadratically perturbed convex functional,
whose solutions exhibit most of the features of the full system we
started with, including the existence of a large class of solutions which
die off weakly at infinity.  In Section~\ref{sec:diff-ineq} we investigate the
asymptotic behavior of solutions to scalar differential equations and
inequalities involving fast oscillating terms.
Section~\ref{sec:oscillating} is devoted to estimates on oscillating
integrals.  Finally, in section~\ref{sec:proof} we put things together
and we conclude the proof of our main results.

\setcounter{equation}{0}
\section{Statements}\label{sec:statements}

Let us consider equation (\ref{eqn:main}) with initial data 
(\ref{eqn:data}). If $A$ is self-adjoint and nonnegative, it is quite 
standard that the problem admits a unique weak global solution
$$u\in C^{1}\left([0,+\infty),H\right)
\cap  C^{0}\left([0,+\infty),D(A^{1/2})\right).$$

Moreover, the classical energy
\begin{equation}
	E(t):=|u'(t)|^{2}+|A^{1/2}u(t)|^{2}
	\label{defn:E}
\end{equation}
is of class $C^{1}$, and its time-derivative satisfies
\begin{equation}
	E'(t)=-2|u'(t)|^{4}
	\quad\quad
	\forall t\geq 0.
	\label{eqn:E'}
\end{equation}

The following is the main result of this paper.

\begin{thm}\label{thm:main}
	Let $H$ be a Hilbert space, and let $A$ be a linear operator on
	$H$ with dense domain $D(A)$.  Let us assume that there exist a
	countable orthonormal basis $\{e_{k}\}$ of $H$ and an increasing sequence
	$\{\lambda_{k}\}$ of positive real numbers such that
	$$Ae_{k}=\lambda_{k}^{2}e_{k}
	\quad\quad
	\forall k\in\n.$$
		
	Let $u(t)$ be the solution to problem
	(\ref{eqn:main})--(\ref{eqn:data}), let $\{u_{0k}\}$ and
	$\{u_{1k}\}$ denote the components of $u_{0}$ and $u_{1}$ with
	respect to the orthonormal basis, and let $\{u_{k}(t)\}$ denote the
	corresponding components of $u(t)$.  Let us consider the set
	\begin{equation}
		J:=\{k\in\n:u_{1k}^{2}+u_{0k}^{2}\neq 0\}.
		\label{defn:J}
	\end{equation}
	
	Then the asymptotic behavior of $u(t)$ and its energy depends on
	$J$ as follows.
	\begin{enumerate}
		\renewcommand{\labelenumi}{(\arabic{enumi})}
		
		\item \emph{(Trivial solution)} If $J=\emptyset$, then
		$u(t)=0$ for every $t\geq 0$ and in particular
		$$\lim_{t\to +\infty}
		t\left(|u'(t)|^{2}+|A^{1/2}u(t)|^{2}\right)=0.$$
	
		\item \emph{(Finite dimensional modes)} If $J$ is a finite set
		with $j$ elements, then $u_{k}(t)=0$ for every $t\geq 0$ and
		every $k\not\in J$.  In addition, for every $k\in J$ there
		exists a real number $\theta_{k,\infty}$ such that
		\begin{equation}
			\lim_{t\to +\infty}\left(
			\sqrt{t}\cdot u_{k}(t)-
			\frac{2}{\sqrt{2j+1}}\cdot\frac{\cos\left(
			\lambda_{k}t+\theta_{k,\infty}\right)}{\lambda_{k}}\,
			\right)=0,
			\label{th:fdm-u}
		\end{equation}
		\begin{equation}
			\lim_{t\to +\infty}\left(
			\sqrt{t}\cdot u_{k}'(t)+
			\frac{2}{\sqrt{2j+1}}\cdot\sin\left(
			\lambda_{k}t+\theta_{k,\infty}\right)\,
			\right)=0,
			\label{th:fdm-u'}
		\end{equation}
		and in particular
		$$\lim_{t\to +\infty}
		t\left(|u'(t)|^{2}+|A^{1/2}u(t)|^{2}\right)=
		\frac{4j}{2j+1}.$$
	
		\item \emph{(Infinite dimensional modes)} If $J$ is
		infinite, then
		$$\lim_{t\to +\infty}
		t\left(|u_{k}'(t)|^{2}+\lambda_{k}^{2}|u_{k}(t)|^{2}\right)=0
		\quad\quad
		\forall k\in\n,$$
		but
		\begin{equation}
			\liminf_{t\to +\infty}
			t\left(|u'(t)|^{2}+|A^{1/2}u(t)|^{2}\right)>0,
			\label{th:liminf-R-main}
		\end{equation}
		and hence $\sqrt{t}\cdot(u(t),u'(t))$ converges to $(0,0)$
		weakly but not strongly.
	\end{enumerate}

\end{thm}

Let us comment on some aspects of Theorem~\ref{thm:main} above.

\begin{rmk}
	\begin{em}
		The result holds true also when $H$ is a finite dimensional 
		Hilbert space, but in that case only the first two options 
		apply.
	\end{em}
\end{rmk}

\begin{rmk}
	\begin{em}
		In the case of finite dimensional modes, let us set
		$$v_{\infty}(t):=\frac{2}{\sqrt{2j+1}}\cdot\sum_{k\in J}
		\frac{\cos\left(\lambda_{k}t+
		\theta_{k,\infty}\right)}{\lambda_{k}}\cdot e_{k}
		\quad\quad
		\forall t\geq 0.$$
		
		It can be verified that $v_{\infty}(t)$ is a solution to the 
		linear homogeneous equation without damping 
		(\ref{eqn:lin-hom}), and that (\ref{th:fdm-u}) and 
		(\ref{th:fdm-u'}) are equivalent to saying that 
		$v_{\infty}(t)$ is the asymptotic profile of $\sqrt{t}\cdot 
		u(t)$, in the sense that
		$$\lim_{t\to +\infty}\left(
		\left|\sqrt{t}\cdot u'(t)-v_{\infty}'(t)\right|^{2}+
		\left|\sqrt{t}\cdot u(t)-v_{\infty}(t)\right|^{2}\right)=0.$$
	\end{em}
\end{rmk}


\begin{rmk}
	\begin{em}
		The assumptions of Theorem~\ref{thm:main} imply in particular that all eigenvalues are simple. Things become more complex if multiplicities are allowed. Let 
		us consider the simplest case where $H$ is a space of 
		dimension~2, and the operator $A$ is the identity.  In this 
		case equation (\ref{eqn:main}) reduces to a system of two 
		ordinary differential equations of the form
		$$\begin{array}{l}
			\ddot{u}+(\dot{u}^{2}+\dot{v}^{2})\dot{u}+u=0,  \\
			\ddot{v}+(\dot{u}^{2}+\dot{v}^{2})\dot{v}+v=0.
		\end{array}$$
		
		If $(v(0),v'(0))=c(u(0),u'(0))$ for some constant $c$, then
		$v(t)=cu(t)$ for every $t\geq 0$, hence there is no
		equipartition of the energy in the limit.  	\end{em}
\end{rmk}

In our second result we consider again the case where $J$ is infinite, and we improve (\ref{th:liminf-R-main}) under a uniform gap condition on eigenvalues (which is satisfied for our model problem (\ref{PDE:nonlocal})).

\begin{thm}\label{thm:energ.} Let $H$, $A$, $\lambda_{k}$, $u(t)$ and $J$ be as in Theorem~\ref{thm:main}. Let us assume in addition that $J$ is infinite and
	\begin{equation}
		\inf_{k\in\n}(\lambda_{k+1}-\lambda_{k})>0.
		\label{hp:diff-lambda}
	\end{equation}

Then it turns out that	
	\begin{equation}
			\lim_{t\to +\infty}
			t\left(|u'(t)|^{2}+|A^{1/2}u(t)|^{2}\right)=2.
			\label{th:lim-R-2-main}
		\end{equation}
\end{thm}

\setcounter{equation}{0}
\section{Basic energy estimates and reduction to ODEs}\label{sec:basic}

In this section we move the first steps in the proof of 
Theorem~\ref{thm:main}. In particular, we prove a basic energy 
estimate and we reduce the problem to a system of countably many 
ordinary differential equations.

\begin{prop}[Basic energy estimate]\label{prop:basic}
	Let $H$, $A$ and $u(t)$ be as in Theorem~\ref{thm:main}. Let us 
	assume that $(u_{0},u_{1})\neq(0,0)$.
	
	Then there exists two positive constants $M_{1}$ and $M_{2}$ such 
	that
	\begin{equation}
		\frac{M_{1}}{1+t}\leq
		|u'(t)|^{2}+|A^{1/2}u(t)|^{2}\leq
		\frac{M_{2}}{1+t}
		\quad\quad
		\forall t\geq 0.
		\label{th:basic-energy}
	\end{equation}
\end{prop}

\paragraph{\textmd{\textit{Proof}}}

Let us consider the classic energy (\ref{defn:E}). From 
(\ref{eqn:E'}) it follows that
$$E'(t)=-2|u'(t)|^{4}\geq -2[E(t)]^{2}
\quad\quad
\forall t\geq 0.$$

Integrating this differential inequality we obtain the estimate from 
below in~(\ref{th:basic-energy}). Since $E'(t)\leq 0$ for every 
$t\geq 0$, we deduce also that
\begin{equation}
	E(t)\leq E(0)
	\quad\quad
	\forall t\geq 0.
	\label{est:EtE0}
\end{equation}

Let us consider now the modified energy
$$F_{\ep}(t):=E(t)+2\ep\langle u(t),u'(t)\rangle E(t),$$
where $\ep$ is a positive parameter. We claim that there exists 
$\ep_{0}>0$ such that 
\begin{equation}
	\frac{1}{2}E(t)\leq F_{\ep}(t)\leq 2E(t)
	\quad\quad
	\forall t\geq 0,\ \forall\ep\in(0,\ep_{0}],
	\label{est:E-F}
\end{equation}
and
\begin{equation}
	F_{\ep}'(t)\leq -\ep[E(t)]^{2}
	\quad\quad
	\forall t\geq 0,\ \forall\ep\in(0,\ep_{0}].
	\label{ineq:F'}
\end{equation}

If we prove these claims, then we set $\ep=\ep_{0}$, and from
(\ref{ineq:F'}) and the estimate from above in (\ref{est:E-F}) we
deduce that
$$F_{\ep_{0}}'(t)\leq -\frac{\ep_{0}}{4}[F_{\ep_{0}}(t)]^{2}
\quad\quad
\forall t\geq 0.$$

An integration of this differential inequality gives that
$$F_{\ep_{0}}(t)\leq\frac{k_{1}}{1+t}
\quad\quad
\forall t\geq 0$$
for a suitable constant $k_{1}$, and hence the estimate from 
below in (\ref{est:E-F}) implies that
$$E(t)\leq 2F_{\ep_{0}}(t)\leq \frac{2k_{1}}{1+t}
\quad\quad
\forall t\geq 0,$$
which proves the estimate from above in~(\ref{th:basic-energy}). 

So we only need to prove (\ref{est:E-F}) and (\ref{ineq:F'}).  The
coerciveness of the operator $A$ implies that
$$\left|2\langle u'(t),u(t)\rangle\right|\leq
|u'(t)|^{2}+|u(t)|^{2}\leq
|u'(t)|^{2}+\frac{1}{\lambda_{1}^{2}}|A^{1/2}u(t)|^{2},$$
so that from (\ref{est:EtE0}) we obtain
\begin{equation}
	\left|2\langle u'(t),u(t)\rangle\right|\leq
	\max\left\{1,\frac{1}{\lambda_{1}^{2}}\right\}E(t)\leq k_{2}
	\quad\quad
	\forall t\geq 0
	\label{est:u-u'}
\end{equation}
for a suitable constant $k_{2}$ depending on initial data.  This
guarantees that (\ref{est:E-F}) holds true when $\ep$ is small enough.

As for (\ref{ineq:F'}), with some computations we obtain that it is 
equivalent to
\begin{eqnarray}
	(2-3\ep)|u'(t)|^{4}+ \ep|A^{1/2}u(t)|^{4}-
	2\ep|u'(t)|^{2}\cdot|A^{1/2}u(t)|^{2} &  & 
	\nonumber  \\
	\noalign{\vspace{1ex}}
	\mbox{}+ 
	6\ep\langle u'(t),u(t)\rangle\cdot|u'(t)|^{4}+ 2\ep\langle
	u'(t),u(t)\rangle\cdot|u'(t)|^{2}\cdot|A^{1/2}u(t)|^{2} & \geq  & 0.
	\label{est:F'-big}
\end{eqnarray}

Keeping (\ref{est:u-u'}) into account, (\ref{est:F'-big}) holds true
if we show that
$$(2-3\ep-3\ep k_{2})|u'(t)|^{4}+ \ep|A^{1/2}u(t)|^{4}-
\ep(2+k_{2})|u'(t)|^{2}\cdot|A^{1/2}u(t)|^{2}\geq 0.$$

The left-hand side is a quadratic form in the variables $|u'(t)|^{2}$ 
and $|A^{1/2}u(t)|^{2}$, and it is nonnegative for all values of the 
variables provided that
$$(2-3\ep-3\ep k_{2})\ep\geq 4\ep^{2}(2+k_{2})^{2},$$
which is clearly true when $\ep$ is small enough. This completes 
the proof.\qed
\bigskip


Proposition~\ref{prop:basic} suggests that $u(t)$ decays as 
$t^{-1/2}$, and  motivates the variable change
$$v(t):=\sqrt{t+1}\cdot u(t)
\quad\quad
\forall t\geq 0.$$

The energy of $v(t)$ is given by
$$|v'(t)|^{2}+|A^{1/2}v(t)|^{2}=
(t+1)|u'(t)|^{2}+
\frac{|u(t)|^{2}}{4(t+1)}+
\langle u'(t),u(t)\rangle+
(t+1)|A^{1/2}u(t)|^{2}.$$

We claim that there exist constants $M_{3}$ and $M_{4}$ such that
\begin{equation}
	0< M_{3}\leq|v'(t)|^{2}+|A^{1/2}v(t)|^{2}\leq M_{4}
	\quad\quad
	\forall t\geq 0.
	\label{est:v-basic}
\end{equation}

The upper estimate being quite clear, we just prove the lower bound. To this end we start by the simple inequality 
\begin{eqnarray*}
|v'(t)|^{2}+|A^{1/2}v(t)|^{2} &  \ge  & (t+1)|u'(t)|^{2}+
\left[\frac{\lambda_1^2}{2}+\frac{1}{4(t+1)}\right]|u(t)|^{2}+
\langle  u'(t), u(t)\rangle  \\ 
&&  \mbox{}+\frac{t+1}{2}|A^{1/2}u(t)|^{2}.
\end{eqnarray*}

On the other hand $$ (t+1)|u'(t)|^{2}+
\left[\frac{\lambda_1^2}{2}+\frac{1}{4(t+1)}\right]|u(t)|^{2}+
\langle  u'(t), u(t)\rangle$$ is obviously greater than or equal to $$ (t+1)|u'(t)|^{2}
+\frac{2\lambda_1^2+ 1}{4(t+1)}|u(t)|^{2}+
\langle  u'(t), u(t)\rangle.$$ 

By decomposing this expression we obtain the inequality
$$ \frac{t+1}{2\lambda_1^2+ 1}|u'(t)|^{2}
+\frac{2\lambda_1^2+ 1}{4(t+1)}|u(t)|^{2}+
\langle  u'(t), u(t)\rangle  + (t+1)\left(1- \frac{1}{2\lambda_1^2+ 1}\right)|u'(t)|^{2}$$ $$ \ge  \frac{2\lambda_1^2}{2\lambda_1^2+ 1}(t+1)|u'(t)|^{2},$$ 
and we end up with 
$$|v'(t)|^{2}+|A^{1/2}v(t)|^{2}\ge \min\left\{\frac{1}{2}, \frac{2\lambda_1^2}{2\lambda_1^2+ 1}\right\}(t+1)(|u'(t)|^{2}+|A^{1/2}u(t)|^{2}),$$ 
which proves the lower bound in (\ref{est:v-basic}) with 
$$ M_3 = \min\left\{\frac{1}{2}, \frac{2\lambda_1^2}{2\lambda_1^2+ 1}\right\} M_1.$$

Starting from (\ref{eqn:main}), with some computations we can verify
that $v(t)$ solves
\begin{equation}
	v''(t)+\left(|v'(t)|^{2}-1\right)\frac{v'(t)}{t+1}
	+Av(t)=g_{1}(t)v(t)+g_{2}(t)v'(t),
	\label{eqn:v}
\end{equation}
where $g_{1}:[0,+\infty)\to\re$ and $g_{2}:[0,+\infty)\to\re$ 
are defined by
$$g_{1}(t):=-\frac{3}{4}\frac{1}{(t+1)^{2}}+
\frac{1}{2}\frac{|v'(t)|^{2}}{(t+1)^{2}}-
\frac{1}{2}\frac{\langle v(t),v'(t)\rangle}{(t+1)^{3}}+
\frac{1}{8}\frac{|v(t)|^{2}}{(t+1)^{4}},$$
$$g_{2}(t):=\frac{\langle v(t),v'(t)\rangle}{(t+1)^{2}}
-\frac{1}{4}\frac{|v(t)|^{2}}{(t+1)^{3}}.$$

Due to (\ref{est:v-basic}), there exists a constant $M_{5}$ such that
\begin{equation}
	|g_{1}(t)|+|g_{2}(t)|\leq\frac{M_{5}}{(t+1)^{2}}
	\quad\quad
	\forall t\geq 0.
	\label{est:g1g2}
\end{equation}

In the sequel we interpret $g_{1}(t)$ and $g_{2}(t)$ as time-dependent
coefficients satisfying this estimate, rather than nonlinear terms.

Let now $\{v_{k}(t)\}$ denote the components of $v(t)$ with respect to 
the orthonormal basis. Then (\ref{eqn:v}) can be rewritten as a system 
of countably many ordinary differential equations of the form
\begin{equation}
	v_{k}''(t)+\left(
	\sum_{i=0}^{\infty}[v_{i}'(t)]^{2}-1\right)\frac{v_{k}'(t)}{t+1}
	+\lambda_{k}^{2}v_{k}(t)=g_{1}(t)v_{k}(t)+g_{2}(t)v_{k}'(t).
	\label{eqn:vk}
\end{equation}

Let us introduce polar coordinates $r_{k}(t)$ and 
$\varphi_{k}(t)$ in such a way that
$$v_{k}(t)=\frac{1}{\lambda_{k}}r_{k}(t)\cos\varphi_{k}(t),
\hspace{3em}
v_{k}'(t)=r_{k}(t)\sin\varphi_{k}(t).$$

In these new variables every second order equation (\ref{eqn:vk}) is 
equivalent to a system of two first order equations of the form (for 
the sake of shortness we do not write explicitly the dependence of 
$r_{k}$ and $\varphi_{k}$ on $t$)
\begin{eqnarray}
	r_{k}' & = & -\left(
	\sum_{i=0}^{\infty}r_{i}^{2}\sin^{2}\varphi_{i}-1\right)
	\frac{r_{k}\sin^{2}\varphi_k}{t+1}+
	\gamma_{k}(t) r_{k}\sin\varphi_k,
	\label{eqn:rk-old}  \\
	\varphi_k' & = & -\lambda_{k}-\left(
	\sum_{i=0}^{\infty}r_{i}^{2}\sin^{2}\varphi_i-1\right)
	\frac{\sin\varphi_k\cos\varphi_k}{t+1}+
	\gamma_{k}(t)\cos\varphi_k,
	\label{eqn:tk-old}
\end{eqnarray}
where 
$$\gamma_{k}(t):=\frac{1}{\lambda_{k}}g_{1}(t)\cos\varphi_{k}(t)+
g_{2}(t)\sin\varphi_{k}(t)
\quad\quad
\forall t\geq 0.$$

In particular, since eigenvalues are bounded from below, from 
(\ref{est:g1g2}) it follows that there exists a constant $M_{6}$ such 
that
\begin{equation}
	|\gamma_{k}(t)|\leq\frac{M_{6}}{(t+1)^{2}}
	\quad\quad
	\forall t\geq 0,\ \forall k\in\n.
	\label{est:gamma-k}
\end{equation}

Finally, we perform one more variable change in order to get rid of
$(t+1)$ in the denominators of equations
(\ref{eqn:rk-old})--(\ref{eqn:tk-old}).  To this end, for every
$k\in\n$ we set
$$\rk(t):=r_{k}(e^{t}-1), 
\hspace{4em}
\tk(t):=\varphi_k(e^{t}-1),$$
and we realize that in these new variables system
(\ref{eqn:rk-old})--(\ref{eqn:tk-old}) reads as
\begin{eqnarray}
	\rk' & = & -\left(
	\sum_{i=0}^{\infty}\ri^{2}\sin^{2}\ti-1\right)
	\rk\sin^{2}\tk+
	\Gamma_{1,k}(t)\rk,
	\label{eqn:rk}  \\
	\tk' & = & -\lambda_{k}e^{t}-\left(
	\sum_{i=0}^{\infty}\ri^{2}\sin^{2}\ti-1\right)
	\sin\tk\cos\tk+
	\Gamma_{2,k}(t),
	\label{eqn:tk}
\end{eqnarray}
where
$$\Gamma_{1,k}(t):=e^{t}\gamma_{k}(e^{t}-1)\sin\theta_{k}(t),
\hspace{3em}
\Gamma_{2,k}(t):=e^{t}\gamma_{k}(e^{t}-1)\cos\theta_{k}(t),$$
so that from (\ref{est:gamma-k}) it follows that there exists a 
constant $M_{7}$ such that
\begin{equation}
	|\Gamma_{1,k}(t)|+|\Gamma_{2,k}(t)|\leq
	M_{7}e^{-t}
	\quad\quad
	\forall t\geq 0,\ \forall k\in\n.
	\label{est:Gamma}
\end{equation}

We observe that $\rk$ can be factored out in the right-hand side of
(\ref{eqn:rk}), and hence either $\rk(t)=0$ for every $t\geq 0$, or
$\rk(t)>0$ for every $t\geq 0$, where the second option applies if and
only if $k$ belongs to the set $J$ defined in (\ref{defn:J}).  We
observe also that the sequence $\rk(t)$ is square-integrable for every
$t\geq 0$, and the square of its norm
\begin{equation}
	R(t):=\sum_{k=0}^{\infty}\rho_{k}^{2}(t)=
	\sum_{k\in J}\rho_{k}^{2}(t) \label{defn:R}
	\end{equation}
satisfies
$$R(t)=\left(|v'(e^{t}-1)|^{2}+
|A^{1/2}v(e^{t}-1)|^{2}\right).$$

In particular, from (\ref{est:v-basic}) it follows that
\begin{equation}
	0<M_{3}\leq R(t)\leq M_{4}
	\quad\quad
	\forall t\geq 0
	\label{est:basic-R}
\end{equation}
for every nontrivial solution.

Finally, we observe that in the new variables Theorem~\ref{thm:main} and Theorem~\ref{thm:energ.} have been reduced to proving the following facts.
\begin{itemize}
	\item  (Finite dimensional modes) If $J$ is a nonempty finite 
	set, then for every $k\in J$ it turns out that
	\begin{equation}
		\lim_{t\to +\infty}\rk(t)=\frac{2}{\sqrt{2j+1}}
		\label{th:lim-rk-finite}
	\end{equation}
	and there
		exists a real number $\theta_{k,\infty}$ such that
	\begin{equation}
		\lim_{t\to +\infty}
		\left(\tk(t)+\lambda_{k}e^{t}\right)= \theta_{k,\infty}.
		\label{th:lim-tk}
	\end{equation}

	\item  (Infinite dimensional modes) If $J$ is infinite, then
	\begin{equation}
		\lim_{t\to +\infty}\rk(t)=0
		\quad\quad
		\forall k\in\n,
		\label{th:lim-rk-infinite}
	\end{equation}
	and under the additional uniform gap assumption (\ref{hp:diff-lambda}) it turns out that 
	\begin{equation}
		\lim_{t\to +\infty}R(t)=2.
		\label{th:lim-R-2}
	\end{equation}
	
	\end{itemize}

\setcounter{equation}{0}
\section{Heuristics}\label{sec:heuristics}

In this section we make some drastic simplifications in equations
(\ref{eqn:rk})--(\ref{eqn:tk}).  These non-rigorous steps lead to a
simplified model, which is then analyzed rigorously in
Theorem~\ref{thm:brutal} below.  The result is that solutions to the
simplified model exhibit all the features stated in
Theorem~\ref{thm:main} and Theorem~\ref{thm:energ.} for solutions to the full system.  Since the
derivation of the simplified model is not rigorous, we can not exploit
Theorem~\ref{thm:brutal} in the study of
(\ref{eqn:rk})--(\ref{eqn:tk}).  Nevertheless, the proof of
Theorem~\ref{thm:brutal} provides a short sketch without
technicalities of the ideas that are involved in the proof of
the main results.

To begin with, in (\ref{eqn:rk}) and (\ref{eqn:tk}) we ignore the
terms with $\Gamma_{1,k}(t)$ and $\Gamma_{2,k}(t)$.  Indeed these
terms are integrable because of (\ref{est:Gamma}), and hence it is
reasonable to expect that they have no influence on the asymptotic
dynamics.  Now let us consider (\ref{eqn:tk}), which seems to suggest
that $\tk(t)\sim - \lambda_{k}e^{t}$.  If this is true, then the
trigonometric terms in (\ref{eqn:rk}) oscillate very quickly, and in
turn this suggests that some homogenization effect takes place.
Therefore, it seems reasonable to replace all those oscillating terms
with their time-averages.

The time-averages can be easily computed to be
\begin{equation}
	\lim_{t\to+\infty}\frac{1}{t}
	\int_{0}^{t}\sin^{2}\left(\lambda e^{s}\right)\,ds=\frac{1}{2}
	\quad\quad
	\forall\lambda>0,
	\label{int:sin2}
\end{equation}
\begin{equation}
	\lim_{t\to+\infty}\frac{1}{t}
	\int_{0}^{t}\sin^{2}\left(\lambda e^{s}\right)\cdot
	\sin^{2}\left(\mu e^{s}\right)\,ds=\frac{1}{4}
	\quad\quad
	\forall\lambda>\mu>0,
	\label{int:sin2-sin2}
\end{equation}
\smallskip
\begin{equation}
	\lim_{t\to+\infty}\frac{1}{t}
	\int_{0}^{t}\sin^{4}\left(\lambda e^{s}\right)\,ds=\frac{3}{8}
	\quad\quad
	\forall\lambda>0.
	\label{int:sin4}
\end{equation}

A comparison of (\ref{int:sin2}) and (\ref{int:sin2-sin2}) reveals
that the two oscillating functions in the integral
(\ref{int:sin2-sin2}) are in some sense independent when
$\lambda\neq\mu$, while (\ref{int:sin4}) shows that this independence 
fails when $\lambda=\mu$. This lack of independence plays a 
fundamental role in the sequel.

After replacing all oscillating coefficients in (\ref{eqn:rk}) with
their time-averages, we are left with the following system of
\emph{autonomous} ordinary differential equations
\begin{equation}
	\rk'=\rk\left(\frac{1}{2}-\frac{3}{8}\rk^{2}
	-\frac{1}{4}\sum_{i\neq k}\ri^{2}\right)=
	\rk\left(\frac{1}{2}-\frac{1}{8}\rk^{2}
	-\frac{1}{4}\sum_{i=0}^{\infty}\ri^{2}\right)
	\label{eqn:rk-avg}
\end{equation}

Quite magically, this system turns out to be the gradient flow of the
functional
\begin{equation}
	\mathcal{F}(\rho):=-\frac{1}{4}\sum_{k=0}^{\infty}\rk^{2}+
	\frac{1}{16}\left(\sum_{k=0}^{\infty}\rk^{2}\right)^{2}+
	\frac{1}{32}\sum_{k=0}^{\infty}\rk^{4},
	\label{defn:functional}
\end{equation}
where $\rho$ belongs to the space $\ell^{2}_{+}$ of square-summable
sequences of nonnegative real numbers.  Since $\mathcal{F}(\rho)$ is a
continuous quadratic perturbation of a convex functional (the sum of the last two
terms), its gradient flow generates a semigroup in $\ell^{2}_{+}$.
Solutions are expected to be asymptotic to stationary points of
$\mathcal{F}(\rho)$.  In addition to the trivial stationary point with
all components equal to~0, all remaining stationary points $\rho$ are
of the form 
$$\rho_{k}:=\left\{
\begin{array}{ll}
	\dfrac{2}{\sqrt{2j+1}}\quad\quad & \mbox{if }k\in J,  \\
	\noalign{\vspace{1ex}}
	0 & \mbox{if }k\not\in J
\end{array}
\right.$$
for some \emph{finite} subset $J\subseteq\n$ with $j$ elements. Incidentally it is not difficult to check that any 
such stationary point is the minimum point of the restriction of 
$\mathcal{F}(\rho)$ to the subset 
\begin{equation}
	W_{J}:=\left\{\rho\in\ell^{2}_{+}:
	\rk=0\mbox{ for every }k\not\in J\right\}.
	\label{defn:WJ}
\end{equation}
Now we show that the asymptotic behavior of solutions to the averaged
system (\ref{eqn:rk-avg}) corresponds to the results announced in our main theorems.

\begin{thm}[Asymptotics for solutions to the homogenized 
	system]\label{thm:brutal}
	Let $\{\rk(t)\}$ be a solution to system (\ref{eqn:rk-avg}) in 
	$\ell^{2}_{+}$, and let $J:=\left\{k\in\n:\rk(0)>0\right\}$.
	
	Then the asymptotic behavior of the solution depends on $J$ as
	follows.
	\begin{enumerate}
		\renewcommand{\labelenumi}{(\arabic{enumi})}
		
		\item  \emph{(Trivial null solution)} If $J=\emptyset$, then 
		$\rho_{k}(t)=0$ for every $k\in\n$ and every $t\geq 0$. 
	
		\item  \emph{(Finite dimensional modes)} If $J$ is a 
		finite set with $j$ elements, then $\rk(t)=0$ for every 
		$k\not\in J$ and every $t\geq 0$, and
		\begin{equation}
			\lim_{t\to +\infty}\rk(t)=\frac{2}{\sqrt{2j+1}}
			\quad\quad
			\forall k\in J.
			\label{th:finite}
		\end{equation}
		
		In other words, in this case the solution leaves in the
		subspace $W_{J}$ defined by (\ref{defn:WJ}), and tends to the
		minimum point of the restriction of $\mathcal{F}(\rho)$ to
		$W_{J}$.
	
		\item  \emph{(Infinite dimensional modes)} If $J$ is 
		infinite, then
		$$\lim_{t\to +\infty}\rk(t)=0
		\quad\quad
		\forall k\in\n,$$
		but
		\begin{equation}
			\lim_{t\to +\infty}\sum_{k=0}^{\infty}\rk^{2}(t)=2,
			\label{th:rho-2}
		\end{equation}
		and in particular the solution tends to 0 weakly but not
		strongly.
	\end{enumerate}
\end{thm}

\paragraph{\textmd{\textit{Proof}}}

First of all we observe that components with null initial datum 
remain null during the evolution, while components with positive 
initial datum remain positive for all subsequent times.

Then we introduce the total energy  $R(t)$ of the solution, defined as in
(\ref{defn:R}).  Moreover, for every pair of indices $h$ and $k$ in
$J$, we consider the ratio
\begin{equation}
	Q_{h,k}(t):=\frac{\rho_{k}(t)}{\rho_{h}(t)}
	\quad\quad
	\forall t\geq 0,
	\label{defn:Q-hk}
\end{equation}
which is well-defined because the denominator never vanishes.

Simple calculations show that
\begin{equation}
	R'(t)=R(t)-\frac{1}{2}R^{2}(t)-\frac{1}{4}
	\sum_{k\in J}\rk^{4}(t)
	\quad\quad
	\forall t\geq 0,
	\label{eqn:R'}
\end{equation}
and 
\begin{equation}
	Q_{h,k}'(t)=\frac{1}{8}\,\rho_{h}^{2}(t)\cdot Q_{h,k}(t)
	\left(1-Q_{h,k}^{2}(t)\right)
	\quad\quad
	\forall t\geq 0.
	\label{eqn:Q'}
\end{equation}

Now we prove some basic estimates on the energy and the quotients, and
then we distinguish the case where all components tend to~0, and the
case where at least one component does not tend to~0.

\subparagraph{\textmd{\textit{Non-optimal energy estimates}}}

We prove that
\begin{equation}
	\frac{4}{3}\leq \liminf_{t\to +\infty}R(t)\leq
	\limsup_{t\to +\infty}R(t)\leq 2.
	\label{ineq:R-trivial}
\end{equation}

Indeed plugging the trivial estimate
$$0\leq\sum_{k\in J}\rk^{4}(t)\leq
\left(\sum_{k\in J}\rk^{2}(t)\right)^{2}$$
into (\ref{eqn:R'}) we obtain that
$$R(t)-\frac{1}{2}R^{2}(t)-\frac{1}{4}R^{2}(t)\leq R'(t) \leq
R(t)-\frac{1}{2}R^{2}(t)
\quad\quad
\forall t\geq 0.$$

Integrating the two differential inequalities we deduce
(\ref{ineq:R-trivial}).

\subparagraph{\textmd{\textit{Uniform boundedness of quotients}}}

We prove that for every $h\in J$ there exists a constant $D_{h}$ such 
that
\begin{equation}
	Q_{h,k}(t)\leq D_{h}
	\quad\quad
	\forall k\in J,\ \forall t\geq 0.
	\label{est:Q-unif}
\end{equation}

We point out that $D_{h}$ is independent of $k$, and actually it can 
be defined as
\begin{equation}
	D_{h}:=\max\left\{1,
	\max\{Q_{h,k}(0):k\in J\}\strut\right\}.
	\label{defn:Mh}
\end{equation}

To this end, it is enough to remark that solutions to (\ref{eqn:Q'})
are decreasing as long as they are greater than~1, and observe that
the inner maximum in (\ref{defn:Mh}) is well defined because for every
fixed $h\in J$ it turns out that $Q_{h,k}(0)\to 0$ as $k\to+\infty$
(because $\rho_{k}(0)\to 0$ as $k\to +\infty$).

\subparagraph{\textmd{\textit{Case where all components vanish in the 
limit}}}

Let us assume that
\begin{equation}
	\lim_{t\to +\infty}\rk(t)=0
	\quad\quad
	\forall k\in J.
	\label{hp:rk->0}
\end{equation}

In this case we prove that $J$ is infinite and (\ref{th:rho-2}) holds 
true.

Let us assume that $J$ is finite.  Then from (\ref{hp:rk->0}) it
follows that $R(t)\to 0$ as $t\to +\infty$, which contradicts the
estimate from below in (\ref{ineq:R-trivial}).

So $J$ is infinite. In order to prove (\ref{th:rho-2}), let us fix 
any index $h_{0}\in J$. From (\ref{est:Q-unif}) we obtain that
$$\sum_{k\in J}\rk^{4}(t)=
\sum_{k\in J}Q_{h_{0},k}^{2}(t)\rho_{h_{0}}^{2}(t)\cdot\rk^{2}(t)\leq
D_{h_{0}}^{2}\cdot\rho_{h_{0}}^{2}(t)\cdot\sum_{k\in J}\rk^{2}(t).$$

Plugging this estimate into (\ref{eqn:R'}) we deduce that
\begin{equation}
	R(t)-\frac{1}{2}R^{2}(t)-
	\frac{1}{4}D_{h_{0}}^{2}\cdot\rho_{h_{0}}^{2}(t)\cdot R(t)\leq R'(t)
	\leq R(t)-\frac{1}{2}R^{2}(t).
	\label{eqn:R'-1}
\end{equation}

Since $\rho_{h_{0}}^{2}(t)\cdot R(t)\to 0$ as $t\to +\infty$, these
two differential inequalities imply (\ref{th:rho-2}) (we refer to
Proposition~\ref{prop:R} below for a more general result).

\subparagraph{\textmd{\textit{Case where at least one component does 
not vanish in the limit}}}

Let us assume that there exists $h_{0}\in J$ such that
\begin{equation}
	\limsup_{t\to +\infty}\rho_{h_{0}}(t)>0.
	\label{hp:rh>0}
\end{equation}

In this case we prove that $J$ is finite and (\ref{th:finite}) holds 
true.

Since $\rho_{h_{0}}(t)$ is Lipschitz continuous (because its
time-derivative is bounded), from (\ref{hp:rh>0}) we deduce that
$$\int_{0}^{+\infty}\rho_{h_{0}}^{2}(t)\,dt=+\infty,$$
and hence from equation (\ref{eqn:Q'}) we conclude that
(we refer to Proposition~\ref{prop:Q1} below for a more general result)
\begin{equation}
	\lim_{t\to +\infty}Q_{h_{0},k}(t)=1
	\quad\quad
	\forall k\in J.
	\label{lim:Q->1}
\end{equation}

We are now ready to prove that $J$ is finite.  Let us assume by
contradiction that this is not the case.  Then for every $n\in\n$
there exists a subset $J_{n}\subseteq J$ with $n$ elements, and hence
$$R(t)\geq
\sum_{k\in J_{n}}\rk^{2}(t)=
\sum_{k\in J_{n}}Q_{h_{0},k}^{2}(t)\rho_{h_{0}}^{2}(t)=
\rho_{h_{0}}^{2}(t)\sum_{k\in J_{n}}Q_{h_{0},k}^{2}(t).$$

When $t\to +\infty$ the last sum tends to $n$ because of 
(\ref{lim:Q->1}), and hence
$$\limsup_{t\to +\infty}R(t)\geq
n\cdot\limsup_{t\to +\infty}\rho_{h_{0}}^{2}(t),$$
which contradicts the estimate from above in (\ref{ineq:R-trivial}) 
when $n$ is large enough. To finish the proof, we now observe that the vector $(\rho_k(t))_{k\in J}$ is a bounded solution of a first order gradient system, so that (cf.\ e.g.\ \cite{SystD}, example 2.2.5 p. 21 or \cite {H-J}, corollary 7.3.1 p. 69) its omega-limit set is made of stationary points only. But the only stationary point satisfying the condition of having all its limiting components positive and equal is the point with all components equal to the right-hand side of (\ref{th:finite}).\qed

\setcounter{equation}{0}
\section{Estimates for differential inequalities}\label{sec:diff-ineq}

In this section we investigate the asymptotic behavior of solutions to
two scalar differential equations, characterized by the presence of
fast oscillating terms. Equations of this form are going to appear 
in the proof of our main results as the equations solved by the energy 
of the solution and by the ratio between two Fourier components.
 \bigskip 

Throughout the text we shall meet oscillatory functions with are not absolutely integrable at infinity 
but have a convergent integral in a weaker sense. 

\begin{defn}[Semi-integrable function] \begin{em}

A function $ f\in C^{0}([t_0, \infty), \re)$ will be called semi-integrable on $[t_0, \infty)$ if the integral 
$$ F(t) : = \int_{t_0}^t f(s)\,ds $$ 
converges to a finite limit as $t$ tends to $+\infty$. In this case the limit will be denoted as $ \int_{t_0}^{+\infty} f(s)\,ds$. 

\end{em}\end{defn}

\begin{rmk}\begin{em}  

A classical example of function which is semi-integrable but not absolutely integrable in $[t_0,+\infty)$ for $t_0>0$  is 
\begin{equation}
f(t)= \frac{\cos(\omega t+ \phi)}{ t^{\alpha}}
\label{ex:semi-int}
\end{equation}
whenever $0<\alpha\leq 1$. Another classical case (Fresnel's integrals) is 
$$f(t)= \cos(\omega t^2+ \phi).$$ 

In the second case the integrability comes from fast oscillations at infinity and the convergence of the integral appears immediately by the change of variable $s= t^2$ which reduces us to (\ref{ex:semi-int}) with $\alpha = 1/2$. The semi-integrable functions that we shall handle are closer to $\cos(ce^{bt})$ in $[0,+\infty)$, in which case the integral can be reduced to (\ref{ex:semi-int}) with $\alpha=1$ by the change of variable $s= e^{bt}$.

\end{em}\end{rmk}

The first equation we consider is actually a differential inequality
which generalizes (\ref{eqn:R'-1}).  It takes the
form
\begin{equation}
	\left|z'(t)-z(t)+\frac{1}{z_{\infty}}\cdot z^{2}(t)-
	\psi_{1}(t)\right|\leq\psi_{2}(t)
	\quad\quad
	\forall t\geq 0.
	\label{hp:R-eqn}
\end{equation}

When $z_{\infty}$ is a positive constant, and
$\psi_{1}(t)\equiv\psi_{2}(t)\equiv 0$, this inequality reduces to an
ordinary differential equation, and it is easy to see that all its
positive solutions tend to $z_{\infty}$ as $t\to +\infty$.  In the
following statement we show that the same conclusion is true under
more general assumption on $\psi_{1}(t)$ and $\psi_{2}(t)$.

\begin{prop}\label{prop:R}
	Let $z_{\infty}$ be a positive constant, and let
	$z:[0,+\infty)\to\re$ be a solution of class $C^{1}$ to the{}
	differential inequality (\ref{hp:R-eqn}). Let us assume that
	\begin{enumerate}
		\renewcommand{\labelenumi}{(\roman{enumi})}
		\item  the function $\psi_{1}:[0,+\infty)\to\re$ is 
		continuous and semi-integrable on $[0,+\infty). $	
		\item  the function $\psi_{2}:[0,+\infty)\to\re$ is 
		continuous and satisfies
		\begin{equation}
			\lim_{t\to +\infty}\psi_{2}(t)=0,
			\label{hp:R-psi2}
		\end{equation}
	
		\item  there exists a constant $c_{0}$ such that
		\begin{equation}
			z(t)\geq c_{0}>0
			\quad\quad
			\forall t\geq 0.
			\label{hp:R-zpos}
		\end{equation}
	\end{enumerate}
	
	Then it turns out that
	\begin{equation}
		\lim_{t\to +\infty}z(t)=z_{\infty}.
		\label{th:R}
	\end{equation}
\end{prop}

\paragraph{\textmd{\textit{Proof}}}

For every $t\geq 0$ let us set
$$x(t):=z(t)-z_{\infty},
\hspace{4em}
a(t):=1+\frac{x(t)}{z_{\infty}}=
\frac{z(t)}{z_{\infty}}.$$

Now (\ref{hp:R-eqn}) is equivalent to the two differential 
inequalities
\begin{equation}
	x'(t)\leq -a(t)x(t)+\psi_{1}(t)+\psi_{2}(t),
	\label{eqn-R-x<}
\end{equation}
\begin{equation}
	x'(t)\geq -a(t)x(t)+\psi_{1}(t)-\psi_{2}(t),
	\label{eqn-R-x>}
\end{equation}
assumption (\ref{hp:R-zpos}) implies that
\begin{equation}
	a(t)\geq\frac{c_{0}}{z_{\infty}}
	\quad\quad
	\forall t\geq 0,
	\label{hp:x-pos}
\end{equation}
and (\ref{th:R}) is equivalent to
\begin{equation}
	\lim_{t\to +\infty}x(t)=0.
	\label{th:R-x}
\end{equation}

Let us set
$$A(t):=\int_{0}^{t}a(\tau)\,d\tau
\quad\quad
\forall t\geq 0,$$
and let us observe that (\ref{hp:x-pos}) implies that $A(t)$ is 
increasing and
\begin{equation}
	\lim_{t\to +\infty}A(t)=+\infty.
	\label{R-lim-A}
\end{equation}

Let us concentrate on the differential inequality (\ref{eqn-R-x<}).
Due to a well-known formula, every solution satisfies
$$x(t)\leq e^{-A(t)}x(0)+
e^{-A(t)}\int_{0}^{t}e^{A(\tau)}\psi_{2}(\tau)\,d\tau+
e^{-A(t)}\int_{0}^{t}e^{A(\tau)}\psi_{1}(\tau)\,d\tau.$$

We claim that the three terms in the right-hand side tend to 0 as
$t\to +\infty$, and hence
\begin{equation}
	\limsup_{t\to +\infty}x(t)\leq 0.
	\label{th:R-limsup}
\end{equation}

This is clear for the first term because of (\ref{R-lim-A}).  Since
$A(t)$ is increasing and tends to $+\infty$, we can apply De
L'H\^{o}pital's rule to the second term.  Keeping (\ref{hp:R-psi2})
and (\ref{hp:x-pos}) into account, we obtain that
$$\lim_{t\to +\infty}\frac{1}{e^{A(t)}}
\int_{0}^{t}e^{A(\tau)}\psi_{2}(\tau)\,d\tau=
\lim_{t\to +\infty}\frac{1}{a(t)e^{A(t)}}
\cdot e^{A(t)}\psi_{2}(t)=0.$$

In order to estimate the third term, let us introduce the function
$$\Psi_{1}(t):=\int_{t}^{+\infty}\psi_{1}(\tau)\,d\tau
\quad\quad
\forall t\geq 0.$$

Due to the semi-integrability of $\psi_{1}(t)$, the function $\Psi_{1}(t)$ is 
well defined and $\Psi_{1}(t)\to 0$ as $t\to +\infty$. Now an 
integration by parts gives that
$$\int_{0}^{t}e^{A(\tau)}\psi_{1}(\tau)\,d\tau=
e^{A(t)}\Psi_{1}(t)-\Psi_{1}(0)-
\int_{0}^{t}a(\tau)e^{A(\tau)}\Psi_{1}(\tau)\,d\tau.$$

The first two terms tend to 0 when multiplied by $e^{-A(t)}$.  As for
the third term, we apply again De L'H\^{o}pital's rule and we conclude
that 
$$\lim_{t\to +\infty}\frac{1}{e^{A(t)}}
\int_{0}^{t}a(\tau)e^{A(\tau)}\Psi_{1}(\tau)\,d\tau= 
\lim_{t\to +\infty}\frac{1}{a(t)e^{A(t)}}
\cdot a(t)e^{A(t)}\Psi_{1}(t)=0.$$

This completes the proof of (\ref{th:R-limsup}).

In an analogous way, from (\ref{eqn-R-x>}) we deduce that
\begin{equation}
	\liminf_{t\to +\infty}x(t)\geq 0.
	\label{th:R-liminf}
\end{equation}

From (\ref{th:R-limsup}) and (\ref{th:R-liminf}) we obtain 
(\ref{th:R-x}), and this completes the proof.\qed
\bigskip

The second equation we consider is a generalization of 
(\ref{eqn:Q'}). It takes the form
\begin{equation}
	z'(t)=\alpha(t)z(t)(1-z^{2}(t))+
	\alpha(t)\beta(t)z^{3}(t)+
	\gamma(t)z(t)
	\quad\quad
	\forall t\geq 0.
	\label{hp:Q1-eqn}
\end{equation}

When $\alpha(t)\equiv 1$ and $\beta(t)\equiv\gamma(t)\equiv 0$, it is 
easy to see that all positive solutions tend to~1 as $t\to +\infty$. 
In the following result we prove the same conclusion under more 
general assumptions on the coefficients.

\begin{prop}\label{prop:Q1}
	Let $z:[0,+\infty)\to(0,+\infty)$ be a positive solution of class $C^{1}$
	to the differential equation (\ref{hp:Q1-eqn}).
	
	Let us assume that
	\begin{enumerate}
		\renewcommand{\labelenumi}{(\roman{enumi})} 
		
		\item the function $\alpha:[0,+\infty)\to(0,+\infty)$ is
		bounded and of class $C^{1}$, and it satisfies
		\begin{equation}
			\int_{0}^{+\infty}\alpha(t)\,dt=+\infty,
			\label{hp:Q1-a}
		\end{equation}
	
		\item  there exists a constant $L_{0}$ such that
		\begin{equation}
			|\alpha'(t)|\leq L_{0}\alpha(t)
			\quad\quad
			\forall t\geq 0,
			\label{hp:Q1-a'}
		\end{equation}
		
		\item the functions $\beta:[0,+\infty)\to\re$ and 
		$\gamma:[0,+\infty)\to\re$ are bounded and semi-integrable.	
	\end{enumerate}
	
	Then it turns out that
	\begin{equation}
		\lim_{t\to +\infty}z(t)=1.
		\label{th:Q1}
	\end{equation}
\end{prop}

\paragraph{\textmd{\textit{Proof}}}

Equation (\ref{hp:Q1-eqn}) is a classical Bernoulli equation, and the
usual variable change $x(t):=[z(t)]^{-2}$ transforms
it into the linear equation
\begin{equation}
	x'(t)=-2(\alpha(t)+\gamma(t))x(t)+2\alpha(t)(1-\beta(t)).
	\label{eqn-Q1-x}
\end{equation}

In the new setting, conclusion (\ref{th:Q1}) is equivalent to proving that
\begin{equation}
	\lim_{t\to +\infty}x(t)=1.
	\label{th:Q1-x}
\end{equation}

In order to avoid plenty of factors~2, with a little abuse of notation
we replace $2\alpha(t)$, $2\beta(t)$, $2\gamma(t)$ with $\alpha(t)$,
$\beta(t)$, $\gamma(t)$.  This does not change the assumptions, but
allows to rewrite (\ref{eqn-Q1-x}) in the simpler form
\begin{equation}
	x'(t)=-(\alpha(t)+\gamma(t))x(t)+\alpha(t)(1-\beta(t)).
	\label{eqn-Q1-x-no2}
\end{equation}

Now we introduce the function
$$A(t):=\int_{0}^{t}\alpha(\tau)\,d\tau
\quad\quad
\forall t\geq 0,$$
and we observe that
\begin{equation}
	\lim_{t\to +\infty}A(t)=+\infty
	\label{Q1-lim-A}
\end{equation}
because of assumption (\ref{hp:Q1-a}). We also introduce the functions
$$B(t):=\int_{t}^{+\infty}\beta(\tau)\,d\tau,
\hspace{3em}
C(t):=\int_{0}^{t}\gamma(\tau)\,d\tau,$$
which are well defined for every $t\geq 0$ as a consequence of assumption 
(iii), and satisfy
\begin{equation}
	\lim_{t\to +\infty}B(t)=0.
	\label{Q1-lim-B}
\end{equation}
\begin{equation}
	\lim_{t\to +\infty}C(t)=:C_{\infty}\in\re.
	\label{Q1-lim-C}
\end{equation}

Every solution to (\ref{eqn-Q1-x-no2}) is given by the well-known formula
\begin{eqnarray}
	x(t) & = & e^{-A(t)-C(t)}x(0)+
	e^{-A(t)-C(t)}\int_{0}^{t}e^{A(\tau)+C(\tau)}\alpha(\tau)\,d\tau
	\nonumber  \\
	 &  & \mbox{}-e^{-A(t)-C(t)}\int_{0}^{t}
	e^{A(\tau)+C(\tau)}\alpha(\tau)\beta(\tau)\,d\tau.
	\nonumber
\end{eqnarray}

We claim that the first and third term tend to 0 as $t\to +\infty$, 
while the second term tends to~1. This would complete the proof of 
(\ref{th:Q1-x}).

The first term tends to 0 because of (\ref{Q1-lim-A}) and
(\ref{Q1-lim-C}).

The second term can be rewritten as
$$e^{-C(t)}\cdot\frac{1}{e^{A(t)}}
\int_{0}^{t}e^{A(\tau)+C(\tau)}\alpha(\tau)\,d\tau.$$

The factor $e^{-C(t)}$ tends to $e^{-C_{\infty}}$.  Since $A(t)$ is
increasing and tends to $+\infty$, we can apply De L'H\^{o}pital's
rule to the second factor.  We obtain that 
$$\lim_{t\to +\infty}\frac{1}{e^{A(t)}}
\int_{0}^{t}e^{A(\tau)+C(\tau)}\alpha(\tau)\,d\tau= \lim_{t\to
+\infty}\frac{1}{\alpha(t)e^{A(t)}}
\cdot e^{A(t)+C(t)}\alpha(t)=e^{C_{\infty}},$$
and this settles the second term.  

In order to compute the limit of the third term, we integrate by
parts.  We obtain that
\begin{eqnarray*}
	\int_{0}^{t}e^{A(\tau)+C(\tau)}\alpha(\tau)\beta(\tau)\,d\tau & = & 
	e^{A(t)+C(t)}\alpha(t)B(t)-\alpha(0)B(0)  \\
	 & & \mbox{}- \int_{0}^{t}e^{A(\tau)+C(\tau)}
	 \left[(\alpha(\tau)+\gamma(\tau))\alpha(\tau)+
	 \alpha'(\tau)\right]B(\tau)\,d\tau.
\end{eqnarray*}

When we multiply by $e^{-A(t)-C(t)}$, the terms in the first line tend
to 0 because of (\ref{Q1-lim-A}) through (\ref{Q1-lim-C}),
and the boundedness of the function $\alpha(t)$.  Thanks to assumption
(\ref{hp:Q1-a'}), the absolute value of the last integral is less than
or equal to
$$\int_{0}^{t}e^{A(\tau)+C(\tau)}
\left(|\alpha(\tau)|+|\gamma(\tau)|+L_{0}\right)
\alpha(\tau)|B(\tau)|\,d\tau.$$

Now we multiply by $e^{-A(t)-C(t)}$, we factor out $e^{-C(t)}$, and we
compute the limit of the rest by exploiting De L'H\^{o}pital's rule as
we did before.  From (\ref{Q1-lim-A}) through (\ref{Q1-lim-C}), and
the boundedness of the functions $\alpha(t)$ and $\gamma(t)$, we
conclude that
$$\hspace{-3em} \lim_{t\to
+\infty}\frac{1}{e^{A(t)}} \int_{0}^{t}e^{A(\tau)+C(\tau)}
\left(|\alpha(\tau)|+|\gamma(\tau)|+L_{0}\right)
\alpha(\tau)|B(\tau)|\,d\tau$$
$$\hspace{3em}=
\lim_{t\to +\infty}
\frac{e^{A(t)+C(t)}\left(|\alpha(t)|+|\gamma(t)|+L_{0}\right)
\alpha(t)|B(t)|}{\alpha(t)e^{A(t)}}=0.$$

This completes the proof of (\ref{th:Q1-x}).\qed
\bigskip

In the third and last result of this section we consider again
equation (\ref{hp:Q1-eqn}).  Let us assume for simplicity that
$\alpha(t)\geq 0$ for every $t\geq 0$, and
$\beta(t)\equiv\gamma(t)\equiv 0$.  These assumptions do not guarantee
that positive solutions tend to~1 as $t\to +\infty$, but nevertheless
they are enough to conclude that all solutions are bounded from above
for $t\geq 0$ (because solutions are decreasing as long as they stay
in the region $z(t)>1$).  In the following result we prove a similar
conclusion under more general assumptions on the coefficients.

\begin{prop}\label{prop:Q2}
	Let $z:[0,+\infty)\to(0,+\infty)$ be a positive solution of class $C^{1}$
	to the differential equation (\ref{hp:Q1-eqn}).
	
	Let us assume that
	\begin{enumerate}
		\renewcommand{\labelenumi}{(\roman{enumi})} 
		
		\item the function
		$\alpha:[0,+\infty)\to(0,+\infty)$ is of class $C^{1}$,
	
		\item the functions $\beta:[0,+\infty)\to\re$ and 
		$\gamma:[0,+\infty)\to\re$ are continuous, 

		\item  there exists a constant $L_{1}$ such that
		\begin{equation}
			\max\left\{\alpha(t),|\alpha'(t)|,
			|\beta(t)|,|\gamma(t)|\right\}\leq L_{1}
			\quad\quad
			\forall t\geq 0,
			\label{hp:Q2-bounded}
		\end{equation}
		
		\item there exists a constant $L_{2}$ such that the following 
		two inequalities
		\begin{equation}
			\left|\int_{t}^{s}\beta(\tau)\,d\tau\right|\leq 
			L_{2}e^{-t},
			\hspace{3em}
			\left|\int_{t}^{s}\gamma(\tau)\,d\tau\right|\leq 
			L_{2}e^{-t}
			\label{hp:Q2-bc}
		\end{equation}
		hold true for every $s\geq t\geq 0$.
	
	\end{enumerate}
	
	Let $t_{0}\geq 0$ be any nonnegative real number such that
	\begin{equation}
		L_{2}\left(1+9L_{1}+32L_{1}^{2}+32L_{1}^{3}\right)
		e^{-t_{0}}<\log 2.
		\label{hp:Q2-t0}
	\end{equation}
	
	Then the following implication holds true
	\begin{equation}
		z(t_{0})\leq 1
		\quad\Longrightarrow\quad
		\sup_{t\geq t_{0}}z(t)\leq 2.
		\nonumber
	\end{equation}
	
\end{prop}

\paragraph{\textmd{\textit{Proof}}}

Let us assume that $z(t_{0})\leq 1$, and let us set
$$t_{2}:=\sup\left\{t\geq t_{0}:
z(\tau)\leq 2\quad\forall\tau\in[t_{0},t]\right\}.$$

If $t_{2}=+\infty$, the result is proved. Let us assume by 
contradiction that this is not the case, and hence $t_{2}<+\infty$. 
Due to the continuity of $z(t)$ and the maximality of $t_{2}$, it 
follows that
\begin{equation}
	z(t_{2})=2.
	\label{Q2:z2}
\end{equation}

Let us set
$$t_{1}:=\inf\left\{t\in[t_{0},t_{2}]: z(\tau)\geq 1\quad
\forall \tau\in[t,t_{2}]\right\}.$$

Then it turns out that $t_{0}\leq t_{1}<t_{2}$, and moreover
\begin{equation}
	z(t_{1})=1
	\label{Q2:z1}
\end{equation}
and
\begin{equation}
	1\leq z(t)\leq 2
	\quad\quad
	\forall t\in[t_{1},t_{2}].
	\label{Q2:z12}
\end{equation}

Due to (\ref{hp:Q2-bounded}) and (\ref{Q2:z12}), from 
(\ref{hp:Q1-eqn}) we deduce that
\begin{equation}
	|z'(t)|\leq 8L_{1}+8L_{1}^{2}
	\quad\quad
	\forall t\in[t_{1},t_{2}].
	\label{Q2:z'}
\end{equation}

Since $z(t)\geq 1$ in $[t_{1},t_{2}]$ and $\alpha(t)$ is positive, 
(\ref{hp:Q1-eqn}) implies also that
$$z'(t)\leq\left(\gamma(t)+\alpha(t)\beta(t)z^{2}(t)\right)z(t)
\quad\quad
\forall t\in[t_{1},t_{2}],$$
which we can integrate as a linear differential inequality. Keeping 
(\ref{Q2:z1}) into account, we find that
$$z(t)\leq\exp\left(\int_{t_{1}}^{t}\gamma(\tau)\,d\tau+
\int_{t_{1}}^{t}\alpha(\tau)\beta(\tau)z^{2}(\tau)\,d\tau\right)
\quad\quad
\forall t\in[t_{1},t_{2}].$$

Now we claim that
\begin{equation}
	\int_{t_{1}}^{t_{2}}\gamma(\tau)\,d\tau+
	\int_{t_{1}}^{t_{2}}\alpha(\tau)\beta(\tau)z^{2}(\tau)\,d\tau
	<\log 2.
	\label{Q2:claim}
\end{equation}

This would imply that $z(t_{2})<2$, thus contradicting (\ref{Q2:z2}).

Due to the second inequality in (\ref{hp:Q2-bc}), we can estimate the
first integral as
\begin{equation}
	\int_{t_{1}}^{t_{2}}\gamma(\tau)\,d\tau\leq
	L_{2}e^{-t_{1}}\leq L_{2}e^{-t_{0}}.
	\label{Q2:claim-1}
\end{equation}

In order to estimate the second integral, we introduce the function
$$B(t):=\int_{t}^{+\infty}\beta(\tau)\,d\tau
\quad\quad
\forall t\geq 0.$$

This function is well defined because of the first inequality 
in (\ref{hp:Q2-bc}), and for the same reason it satisfies
\begin{equation}
	B(t)\leq L_{2}e^{-t}
	\quad\quad
	\forall t\geq 0.
	\label{Q2:est-B}
\end{equation}

Now an integration by parts gives that
\begin{eqnarray*}
	\int_{t_{1}}^{t_{2}}\alpha(\tau)\beta(\tau)z^{2}(\tau)\,d\tau & = & 
	\alpha(t_{2})z^{2}(t_{2})B(t_{2})-\alpha(t_{1})z^{2}(t_{1})B(t_{1}) \\
	 & & -\int_{t_{1}}^{t_{2}}B(\tau)\left(
	 \alpha'(\tau)z^{2}(\tau)+2\alpha(\tau)z(\tau)z'(\tau)\right)\,d\tau.
\end{eqnarray*}

From (\ref{hp:Q2-bounded}), (\ref{Q2:z2}), (\ref{Q2:z1}), and
(\ref{Q2:est-B}) it follows that
\begin{eqnarray}
	\left|\alpha(t_{2})z^{2}(t_{2})B(t_{2})-
	\alpha(t_{1})z^{2}(t_{1})B(t_{1})\right| & \leq & 
	L_{1}\cdot 4\cdot L_{2}e^{-t_{2}}+
	L_{1}\cdot 1\cdot L_{2}e^{-t_{1}}
	\nonumber  \\
	\noalign{\vspace{0.5ex}}
	 & \leq & 5L_{1}L_{2}e^{-t_{0}}.
	\label{Q2:est-1}
\end{eqnarray}

From (\ref{hp:Q2-bounded}), (\ref{Q2:z12}), (\ref{Q2:z'}) and
(\ref{Q2:est-B}) it follows that
\begin{eqnarray}
	\left|B(\tau)\left(\alpha'(\tau)z^{2}(\tau)+
	2\alpha(\tau)z(\tau)z'(\tau)\right)\right| & \leq & 
	L_{2}e^{-\tau}\left(4L_{1}+32L_{1}(L_{1}^{2}+L_{1})\right)
	\nonumber  \\
	\noalign{\vspace{0.5ex}}
	 & \leq & 4L_{2}\left(L_{1}+8L_{1}^{2}+8L_{1}^{3}\right)e^{-\tau}
	\label{Q2:est-2}
\end{eqnarray}
for every $\tau\in[t_{1},t_{2}]$. From (\ref{Q2:est-1}) and 
(\ref{Q2:est-2}) it follows that
\begin{equation}
	\int_{t_{1}}^{t_{2}}\alpha(\tau)\beta(\tau)z^{2}(\tau)\,d\tau\leq
	L_{2}\left(9L_{1}+32L_{1}^{2}+32L_{1}^{3}\right)e^{-t_{0}}.
	\label{Q2:claim-2}
\end{equation}

Adding (\ref{Q2:claim-1}) and (\ref{Q2:claim-2}), and taking
assumption (\ref{hp:Q2-t0}) into account, we obtain (\ref{Q2:claim}).
This completes the proof.\qed

\setcounter{equation}{0}
\section{Estimates on oscillating integrals}\label{sec:oscillating}

In the three results of this section we prove the convergence of
some oscillating integrals and series of oscillating integrals.  We need
these estimates in the proof of our main result when we deal with the
trigonometric terms of (\ref{eqn:rk}) and (\ref{eqn:tk}).

\begin{lemma}\label{lemma:osc-int}
	Let $\alpha>0$, let $L_{3}\geq 0$, and let $\psi:[0,+\infty)\to\re$ 
	be a function of class $C^{1}$ such that
	\begin{equation}
		|\psi'(t)|\leq L_{3}
		\quad\quad
		\forall t\geq 0.
		\label{hp:osc-lip}
	\end{equation}
	
	Then for every $s\geq t\geq 0$ it turns out that
	\begin{equation}
		\left|\int_{t}^{s}
		\cos(\alpha e^{\tau}+\psi(\tau))\,d\tau\right|\leq
		\frac{3+L_3}{\alpha e^{t}}.
		\label{th:osc}
	\end{equation}
\end{lemma}

\paragraph{\textmd{\textit{Proof}}}

We introduce the complex valued functions 
$$g(t):=\exp\left(i \alpha e^{t}\right),
\hspace{4em} 
f(t) :=  \exp\left(i \psi(t)\right),$$ 
so that clearly 
$$ \left|\int_{t}^{s}\cos(\alpha e^{\tau}+\psi(\tau))\,d\tau\right|\leq \left|\int_{t}^{s}\exp(i[\alpha e^{\tau}+\psi(\tau)])\,d\tau\right| = \left|\int_{t}^{s}g(\tau)f(\tau)\,d\tau\right|.$$ 

Now we have 
$$\int_{t}^{s}g(\tau)f(\tau)\,d\tau = \int_{t}^{s}g'(\tau) \frac{1}{i\alpha}e^{-\tau}f(\tau)\,d\tau $$ 
$$= \frac{1}{i\alpha}\left[g(s)f(s)e^{-s}- g(t)f(t)e^{-t} -
\int_{t}^{s}g(\tau)(f'(\tau) -f(\tau)) e^{-\tau}d\tau\right],$$ 
yielding the immediate estimate
$$ \left|\int_{t}^{s}g(\tau)f(\tau)\,d\tau\right| \leq \frac{3 + L_3}{\alpha}e^{-t},$$
which implies (\ref{th:osc})\qed
\bigskip

Lemma~\ref{lemma:osc-int} can also be viewed as a special case of the following result.

\begin{lemma}\label{lemma:osc-prod}
	Let $g:[0,+\infty)\to\mathbb{C}$ be a continuous function, and let $f:[0,+\infty)\to\mathbb{C}$ be a function of class $C^{1}$. 
	Let us assume that there exist two constants $L_{4}$ and $L_{5}$ such that
	\begin{equation}
		\left|\int_{t}^{s}g(\tau)\,d\tau\right|\leq
		L_{4}e^{-t}
		\quad\quad
		\forall s\geq t\geq 0,
		\label{hp:OP-a}
	\end{equation}
	\begin{equation}
		\max\left\{|f(t)|,|f'(t)|\right\}\leq L_{5}
		\quad\quad
		\forall t\geq 0.
		\nonumber
	\end{equation}
	
	Then it turns out that
	\begin{equation}
		\left|\int_{t}^{s}g(\tau)f(\tau)\,d\tau\right|\leq
		3L_{4}L_{5}e^{-t}
		\quad\quad
		\forall s\geq t\geq 0.
		\label{th:osc-prod}
	\end{equation}
\end{lemma}

\paragraph{\textmd{\textit{Proof}}}

Let us introduce the function
$$G(t):=\int_{t}^{+\infty}g(\tau)\,d\tau
\quad\quad
\forall t\geq 0.$$

This function is well defined because of assumption (\ref{hp:OP-a}), 
and for the same reason it satisfies
\begin{equation}
	|G(t)|\leq L_{4}e^{-t}
	\quad\quad
	\forall t\geq 0.
	\nonumber
\end{equation}

Integrating by parts the left-hand side of (\ref{th:osc-prod}) we find
that
$$\int_{t}^{s}g(\tau)f(\tau)\,d\tau=
G(s)f(s)-G(t)f(t)-\int_{t}^{s}G(\tau)f'(\tau)\,d\tau.$$

At this point our assumptions imply that
\begin{eqnarray*}
	\left|\int_{t}^{s}g(\tau)f(\tau)\,d\tau\right| & \leq & 
	|G(s)|\cdot|f(s)|+|G(t)|\cdot|f(t)|+
	\int_{t}^{s}|G(\tau)|\cdot|f'(\tau)|\,d\tau\\
	 & \leq & L_{4}e^{-s}\cdot L_{5}+
	 L_{4}e^{-t}\cdot L_{5}+ 
	 \int_{t}^{s}L_{4}e^{-\tau}\cdot L_{5}\,d\tau   \\
	 & \leq & 3L_{4}L_{5}e^{-t},
\end{eqnarray*}
which proves (\ref{th:osc-prod}).\qed
\bigskip 

The next lemma extends the previous estimates to some series of functions.
 
\begin{lemma}\label{lemma:osc-series}
	Let $g_{k}:[0,+\infty)\to\re$ be a sequence of continuous
	functions, and let $f_{k}:[0,+\infty)\to\re$ be a sequence of
	functions of class $C^{1}$.
	
	Let us assume that  the two series of functions
	$$\sum_{k=0}^{\infty}f_{k}(t),
	\hspace{3em}
	\sum_{k=0}^{\infty}f_{k}'(t)$$
	are normally convergent on compact subsets of $[0,+\infty)$, and
	that there exist three constants $L_{6}$, $L_{7}$, and $L_{8}$ such
	that
	\begin{equation}
		|g_{k}(t)|\leq L_{6}
		\quad\quad
		\forall t\geq 0,\quad\forall k\in\n,
		\label{hp:OS-a-bd}
	\end{equation}
	\begin{equation}
		\left|\int_{t}^{s}g_{k}(\tau)\,d\tau\right|\leq
		L_{7}e^{-t}
		\quad\quad
		\forall s\geq t\geq 0,\quad\forall k\in\n,
		\label{hp:OS-a}
	\end{equation}
	and
	\begin{equation}
		\max\left\{\sum_{k=0}^{\infty}|f_{k}(t)|,
		\sum_{k=0}^{\infty}|f_{k}'(t)|\right\}\leq L_{8}
		\quad\quad
		\forall t\geq 0.
		\label{hp:OS-f}
	\end{equation}
	
	Then the series
	\begin{equation}
		\sum_{k=0}^{\infty}g_{k}(t)f_{k}(t)
		\label{def:series-af}
	\end{equation}
	is normally convergent on compact subsets of $[0,+\infty)$, and 
	its sum satisfies
	\begin{equation}
		\left|\int_{t}^{s}
		\left(\sum_{k=0}^{\infty}g_{k}(\tau)f_{k}(\tau)\right)
		\,d\tau\right|\leq
		3L_{7}L_{8}e^{-t}
		\quad\quad
		\forall s\geq t\geq 0.
		\label{th:osc-series}
	\end{equation}
\end{lemma}

\paragraph{\textmd{\textit{Proof}}}

In analogy with the proof of Lemma~\ref{lemma:osc-prod}, we introduce
the functions 
$$G_{k}(t):=\int_{t}^{+\infty}g_{k}(\tau)\,d\tau.$$

We observe that they are well defined because of assumption
(\ref{hp:OS-a}), and they satisfy
\begin{equation}
	|G_{k}(t)|\leq L_{7}e^{-t}
	\quad\quad
	\forall t\geq 0,\ \forall k\in\n.
	\label{est:Gk}
\end{equation}

From assumption (\ref{hp:OS-a-bd}) it follows that 
$$\sup_{t\in K}|g_{k}(t)f_{k}(t)|\leq
L_{6}\sup_{t\in K}|f_{k}(t)|
\quad\quad
\forall k\in\n$$
for every compact set $K\subseteq[0,+\infty)$.  As a consequence, the
normal convergence in $K$ of the series (\ref{def:series-af}) follows
from the normal convergence in $K$ of the series with general term
$f_{k}(t)$.  Due to normal convergence, we can exchange series and
integrals in the left-hand side of (\ref{th:osc-series}), and deduce
that
$$\left|\int_{t}^{s}
\left(\sum_{k=0}^{\infty}g_{k}(\tau)f_{k}(\tau)\right)
\,d\tau\right|=
\left|\sum_{k=0}^{\infty}\int_{t}^{s}
g_{k}(\tau)f_{k}(\tau)\,d\tau\right|\leq
\sum_{k=0}^{\infty}\left|\int_{t}^{s}
g_{k}(\tau)f_{k}(\tau)\,d\tau\right|.$$

Now we integrate by parts each term of the series, and we exploit
(\ref{est:Gk}) in analogy with what as we did before in the proof of
Lemma~\ref{lemma:osc-prod}.  We obtain that
$$\left|\int_{t}^{s}g_{k}(\tau)f_{k}(\tau)\,d\tau\right|\leq
L_{7}e^{-s}|f_{k}(s)|+L_{7}e^{-t}|f_{k}(t)|+
L_{7}\int_{t}^{s}e^{-\tau}|f_{k}'(\tau)|\,d\tau$$
for every $k\in\n$. When we sum over $k$, from (\ref{hp:OS-f}) we 
deduce that
\begin{equation}
	\sum_{k=0}^{\infty}L_{7}e^{-s}|f_{k}(s)|=
	L_{7}e^{-s}\sum_{k=0}^{\infty}|f_{k}(s)|\leq L_{7}L_{8}e^{-t},
	\label{OS-1}
\end{equation}
and analogously
\begin{equation}
	\sum_{k=0}^{\infty}L_{7}e^{-t}|f_{k}(t)|\leq L_{7}L_{8}e^{-t}.
	\label{OS-2}
\end{equation}

As for the sum of integrals, we first observe that the normal 
convergence on compact subsets of $[0,+\infty)$ of the series with 
general term $f_{k}'(t)$ implies an analogous convergence of the series
$$\sum_{k=0}^{\infty}e^{-\tau}|f_{k}'(\tau)|.$$

Therefore, we can exchange once again series and integrals.  Keeping
(\ref{hp:OS-f}) into account, this leads to
\begin{eqnarray}
	\sum_{k=0}^{\infty}L_{7}\int_{t}^{s}e^{-\tau}|f_{k}'(\tau)|\,d\tau & = & 
	L_{7}\int_{t}^{s}\left(
	\sum_{k=0}^{\infty}e^{-\tau}|f_{k}'(\tau)|\right)\,d\tau 
	\nonumber  \\
	 & \leq & L_{7}\int_{t}^{s}L_{8}e^{-\tau}\,d\tau
	 \nonumber  \\
	 & \leq & L_{7}L_{8}\,e^{-t}.
	 \label{OS-3}
\end{eqnarray}

At this point, (\ref{th:osc-series}) follows from (\ref{OS-1}),
(\ref{OS-2}), and (\ref{OS-3}).\qed

\setcounter{equation}{0}
\section{Proof of the main results}\label{sec:proof}

\subsection*{Equations for the energy and quotients}

\paragraph{\textmd{\textit{Preliminary estimates on components}}}

Let us consider the notations introduced in section~\ref{sec:basic},
where we reduced ourselves to proving (\ref{th:lim-rk-finite}) through
(\ref{th:lim-R-2}).  In this first paragraph we derive some
$k$-independent estimates on $\rk(t)$ and $\tk(t)$ that are needed
several times in the sequel. The constants $M_{8}$, \ldots, 
$M_{23}$ we introduce hereafter depend on the solution (as 
the constants $M_{1}$, \ldots, $M_{7}$ of section~\ref{sec:basic}), 
but they do not depend on $k$. First of all, from (\ref{defn:R}) and  (\ref{est:basic-R}) it follows  that 
$$\sum_{k=0}^{\infty}\rho_{k}^{2}(t)\leq M_{4}$$ and in particular we find
\begin{equation}
	\rho_{k}(t)\leq M_{8}
	\quad\quad
	\forall t\geq 0,\quad\forall k\in\n,
	\label{est:rk}
\end{equation}
and
\begin{equation}
	\sum_{k=0}^{\infty}\rho_{k}^{2}(t)\sin^{2}\theta_{k}(t)\leq M_{4}.
	\label{est:R-sin}
\end{equation}

From this estimate and (\ref{est:Gamma}) it follows that
\begin{equation}
	|\rho_{k}'(t)|\leq M_{9}\rho_{k}(t)
	\quad\quad
	\forall t\geq 0,\quad\forall k\in\n.
	\label{est:rk'-rk}
\end{equation}

This implies in particular that
\begin{equation}
	\sum_{k=0}^{\infty}\left[\rho_{k}'(t)\right]^{2}\leq M_{10}
	\quad\quad
	\forall t\geq 0,
	\label{est:rk'-series}
\end{equation}
and
\begin{equation}
	|\rho_{k}'(t)|\leq M_{11}
	\quad\quad
	\forall t\geq 0,\quad\forall k\in\n.
	\label{est:rk'}
\end{equation}

Moreover, from (\ref{est:rk'-rk}) it follows that
\begin{equation}
	\rho_{k}(t)\leq\rho_{k}(0)e^{M_{9}t}
	\quad\quad
	\forall t\geq 0,\ \forall k\in\n.
	\label{est:rk-rk0}
\end{equation}

Let us consider now the series
\begin{equation}
	\sum_{k=0}^{\infty}\rho_{k}^{m}(t),
	\hspace{4em}
	\sum_{k=0}^{\infty}[\rho_{k}^{m}(t)]',
	\nonumber
\end{equation}
where $m\geq 2$ is a fixed exponent (in the sequel we need only the
cases $m=2$ and $m=4$).  From the previous estimates it follows that
\begin{equation}
	\sum_{k=0}^{\infty}\rho_{k}^{m}(t)\leq M_{12},
	\hspace{4em}
	\sum_{k=0}^{\infty}\left|[\rho_{k}^{m}(t)]'\strut\right|\leq M_{12},
	\label{est:series-m}
\end{equation}
where of course the constant $M_{12}$ depends also on $m$.  Moreover,
from (\ref{est:rk-rk0}) and the square-integrability of the sequence
$\rho_{k}(0)$, it follows that both series are normally convergent on 
compact subsets of $[0,+\infty)$.

We stress that we can not hope that these series are normally 
convergent in $[0,+\infty)$, even when $m=2$. Indeed normal 
convergence would imply uniform convergence, and hence the 
possibility to exchange the series and the limit as $t\to +\infty$, 
while the conclusion of Theorem~\ref{thm:main} says that this is not 
the case, at least when $J$ is an infinite set.

Finally, plugging (\ref{est:Gamma}) and (\ref{est:R-sin}) into
(\ref{eqn:tk}), after integration we obtain that
\begin{equation}
	\tk(t)=-\lambda_{k}e^{t}-\psi_{k}(t)
	\label{defn:psik}
\end{equation}
for a suitable function $\psi_{k}:[0,+\infty)\to\re$ of class $C^{1}$
satisfying
\begin{equation}
	|\psi_{k}'(t)|\leq M_{13}
	\quad\quad
	\forall t\geq 0,\ \forall k\in\n.
	\label{est:psik'}
\end{equation}

\paragraph{\textmd{\textit{Estimates on trigonometric coefficients}}}

For every $k\in\n$ we set
\begin{equation}
	a_{k}(t):=\sin^{2}\theta_{k}(t)-\frac{1}{2},
	\hspace{3em}
	b_{k}(t):=\sin^{4}\theta_{k}(t)-\frac{3}{8},
	\nonumber
\end{equation}
and for every $k\neq h$ we set
\begin{equation}
	c_{h,k}(t):=\sin^{2}\theta_{h}(t)\sin^{2}\theta_{k}(t)
	-\frac{1}{4}.
	\nonumber
\end{equation}

These functions represent the corrections we have to take into 
account when we approximate the trigonometric functions with their 
time-average, as we did at the beginning of 
section~\ref{sec:heuristics}.

It is easy to see that 
\begin{equation}
	\sup\{|a_{k}(t)|,|b_{k}(t)|,|c_{h,k}(t)|\}\leq 1
	\quad\quad
	\forall t\geq 0,
	\label{est:abc-bd}
\end{equation}
where the supremum is taken over all admissible indices or pairs of
indices.  Now we claim that
\begin{equation}
	\left|\int_{t}^{s}a_{k}(\tau)\,d\tau\right|\leq M_{14}e^{-t}
	\quad\quad
	\forall s\geq t\geq 0,\ \forall k\in\n,
	\label{est:main-a}
\end{equation}
\begin{equation}
	\left|\int_{t}^{s}b_{k}(\tau)\,d\tau\right|\leq M_{15}e^{-t}
	\quad\quad
	\forall s\geq t\geq 0,\ \forall k\in\n,
	\label{est:main-b}
\end{equation}
and 
\begin{equation}
\left|\int_{t}^{s}c_{h,k}(\tau)\,d\tau\right|\leq 
M_{16}\left(1+\frac{1}{|\lambda_{k}-\lambda_{h}|}\right)e^{-t}
	\quad\quad
	\forall s\geq t\geq 0,\ \forall h\neq k.
	\label{est:main-c}\end{equation}

In order to prove (\ref{est:main-a}), we just observe that
$$a_{k}(t)=-\frac{1}{2}\cos(2\tk(t)),$$
and hence by (\ref{defn:psik})
$$a_{k}(t)=-\frac{1}{2}
\cos\left(-2\lambda_{k}e^{t}-2\psi_{k}(t)\right)=-\frac{1}{2}
\cos\left(2\lambda_{k}e^{t}+2\psi_{k}(t)\right).$$

Thanks to (\ref{est:psik'}), the assumptions of 
Lemma~\ref{lemma:osc-int} are satisfied with $\alpha:=2\lambda_{k}$, 
$L_{3}:=2M_{13}$, 
and $\psi(t):=\psi_{k}(t)$. Thus we obtain that
$$\left|\int_{t}^{s}a_{k}(\tau)\,d\tau\right|\leq\frac{3+2M_{13}}{2\lambda_{k}}e^{-t}\leq M_{17}e^{-t},$$
where in the last inequality we exploited that all eigenvalues are
larger than a fixed positive constant.

The proof of (\ref{est:main-b}) is analogous, just starting from the 
trigonometric identity
$$b_{k}(t)=-\frac{1}{2}\cos(2\tk(t))+\frac{1}{8}\cos(4\tk(t)).$$

Also the proof of (\ref{est:main-c}) is analogous, but in this case the 
trigonometric identity is
$$c_{h,k}=-\frac{1}{4}\cos(2\theta_{h})-
\frac{1}{4}\cos(2\theta_{k})
+\frac{1}{8}\cos(2\theta_{h}+2\theta_{k})
+\frac{1}{8}\cos(2\theta_{h}-2\theta_{k}).$$

All the four terms can be treated through Lemma~\ref{lemma:osc-int}, 
but now in the last term the differences between eigenvalues are involved. As a consequence, for the last term we obtain an estimate of the form
$$\left|\int_{t}^{s}\cos\left(
2\theta_{h}(\tau)-2\theta_{k}(\tau)\right)\,d\tau\right|\leq
\frac{3+ 4M_{13}}{2|\lambda_{k}-\lambda_{h}|}e^{-t}.$$

If we want this estimate to be uniform for $k\neq h$, we have to 
assume that differences between eigenvalues are bounded away from~0, 
and this is exactly the point where assumption (\ref{hp:diff-lambda}) 
comes into play in the proof of Theorem \ref{thm:energ.}.

\paragraph{\textmd{\textit{Equation for the energy}}}

Let $R(t)$ be the total energy as defined in  (\ref{defn:R}).  We claim that $R(t)$ solves a
differential equation of the form
\begin{equation}
	R'(t)=R(t)-\frac{1}{2}R^{2}(t)-
	\frac{1}{4}\sum_{k=0}^{\infty}\rho_{k}^{4}(t)+\mu_1(t)+\mu_2(t),
	\label{eqn:R'-full}
\end{equation}
where  (for the sake of 
shortness, we do not write the explicit dependence on $t$ in the 
right-hand sides)
\begin{equation}
	\mu_1(t):=2\sum_{k=0}^{\infty}\left(
	\Gamma_{1,k}\rho_{k}^{2}+
	a_{k}\rho_{k}^{2}-b_{k}\rho_{k}^{4}\right),
\quad\quad	
\mu_2(t):=-2\sum_{k=0}^{\infty}\left(\rho_{k}^{2}\sum_{i\neq k}c_{i,k}\rho_{i}^{2}\right).
	\label{defn:f-g}
\end{equation}

We also claim that $\mu_1(t)$ satisfies
\begin{equation}
	\left|\int_{t}^{s}\mu_{1}(\tau)\,d\tau\right|\leq M_{18}e^{-t}
	\quad\quad
	\forall s\geq t\geq 0.
	\label{est:osc-R}
\end{equation}

The verification of (\ref{eqn:R'-full}) is a lengthy but elementary 
calculation, which starts by writing
$$R'(t)=2\sum_{k=0}^{\infty}\rho_{k}(t)\rho_{k}'(t)$$
and by replacing $\rho_{k}'(t)$ with the right-hand side of 
(\ref{eqn:rk}). The crucial point is that when computing the product
$$\rho_{k}^{2}\sin^{2}\theta_{k}\cdot
\sum_{i=0}^{\infty}\rho_{i}^{2}\sin^{2}\theta_{i},$$
one has to isolate the term of the series with $i=k$. In this way the 
product becomes
$$\rho_{k}^{4}\sin^{4}\theta_{k}+\rho_{k}^{2}
\sum_{i\neq k}\rho_{i}^{2}\sin^{2}\theta_{i}\sin^{2}\theta_{k},$$
and now one can express $\sin^{4}\theta_{k}$ in terms of $b_{k}$, and 
$\sin^{2}\theta_{i}\sin^{2}\theta_{k}$ in terms of $c_{i.k}$. The 
rest is straightforward algebra.

The proof of (\ref{est:osc-R}) follows from several applications 
of Lemma~\ref{lemma:osc-series} with different choices of $f_{k}(t)$ 
and $g_{k}(t)$.
\begin{itemize}
	\item  In the case of the term $\Gamma_{1,k}\rho_{k}^{2}$ we 
	choose $f_{k}(t):=\rk^{2}(t)$ and $g_{k}(t):=\Gamma_{1,k}(t)$. 
	Indeed the assumptions on $f_{k}(t)$ follow from 
	(\ref{est:series-m}) with $m=2$ and from the normal convergence of the 
	same series on compact subsets of $[0,+\infty)$, while the 
	assumptions on $g_{k}(t)$ follows from (\ref{est:Gamma}).

	\item In the case of the term $a_{k}\rho_{k}^{2}$ we choose
	$f_{k}(t):=\rk^{2}(t)$ and $g_{k}(t):=a_{k}(t)$.  The assumptions on
	$f_{k}(t)$ are satisfied as before, while those on $g_{k}(t)$
	follow from (\ref{est:abc-bd}) and (\ref{est:main-a}).

	\item In the case of the term $b_{k}\rho_{k}^{4}$ we choose
	$f_{k}(t):=\rk^{4}(t)$ and $g_{k}(t):=b_{k}(t)$.  Now we need the
	estimates for the series (\ref{est:series-m}) with $m=4$ in order
	to verify the assumptions on $f_{k}(t)$, and (\ref{est:abc-bd})
	and (\ref{est:main-b}) in order to provide the requires estimates
	on $g_{k}(t)$.


\end{itemize}

\paragraph{\textmd{\textit{Equation for quotients}}}

For every pair of indices $h$ and $k$ in $J$ we consider the ratio
$Q_{h,k}(t)$ introduced in (\ref{defn:Q-hk}).  We remind that
components with indices in $J$ never vanish, and therefore the
quotient is well defined and positive for every $t\geq 0$. After some 
lengthy calculations we obtain \begin{equation}
	Q_{h,k}'(t)=\alpha_{h}(t)Q_{h,k}(t)\left(1-Q_{h,k}^{2}(t)\right)
	+\alpha_{h}(t)\beta_{h,k}(t)Q_{h,k}^{3}(t)+
	\gamma_{h,k}(t)Q_{h,k}(t),
	\label{eqn:Q'-true}
\end{equation}
where
\begin{equation}
	\alpha_{h}(t):=\frac{1}{8}\rho_{h}^{2}(t),
	\nonumber
\end{equation}
\begin{equation}
	\beta_{h,k}(t)=8(c_{h,k}(t)-b_{k}(t)),
	\nonumber
\end{equation}
\begin{equation}
	\gamma_{h,k}(t):=a_{k}-a_{h}+\Gamma_{1,k}-\Gamma_{1,h}
	+\rho_{h}^{2}(b_{h}-c_{h,k})
	+\sum_{i\not\in\{h,k\}}\rho_{i}^{2}(c_{i,h}-c_{i,k}).
	\nonumber
\end{equation}

We observe that the first term of equation (\ref{eqn:Q'-true}) is the
same as in equation (\ref{eqn:Q'}), which was derived by neglecting
all the rest.

We claim that 
\begin{equation}
	\sup\left\{|\alpha_{h}(t)|,|\alpha_{h}'(t)|,
	|\beta_{h,k}(t)|,|\gamma_{h,k}(t)|\strut\right\}\leq M_{19}
	\quad\quad
	\forall t\geq 0,
	\label{est:sup}
\end{equation}
where the supremum is taken over all admissible indices or pairs of
indices, and that
\begin{equation}
	\left|\int_{t}^{s}\beta_{h,k}(\tau)\,d\tau\right|
	\leq M_{20}\left(1+\frac{1}{|\lambda_{k}-\lambda_{h}|}\right)e^{-t},
	\label{est:int-beta}
\end{equation}
\begin{equation}
	\left|\int_{t}^{s}\gamma_{h,k}(\tau)\,d\tau\right|
	\leq M_{21}\left(1+\frac{1}{|\lambda_{k}-\lambda_{h}|}+
	\sup_{i\not\in\{h,k\}}\left(\frac{1}{|\lambda_{i}-\lambda_{k}|}+\frac{1}{|\lambda_{i}-\lambda_{h}|}\right)\right)e^{-t}
	\label{est:int-gamma}
\end{equation}
for every pair of admissible indices and every $s>t\geq 0$. We point out that the supremum in (\ref{est:int-gamma}) is finite because the sequence of eigenvalues is increasing.

Estimate (\ref{est:sup}) follows from (\ref{est:rk}) and
(\ref{est:rk'}) in the case of $\alpha_{h}(t)$ and $\alpha_{h}'(t)$,
from (\ref{est:abc-bd}) in the case of $\beta_{h,k}(t)$, and from
(\ref{est:abc-bd}), (\ref{est:Gamma}) and (\ref{est:basic-R}) in the
case of $\gamma_{h,k}(t)$.  

Estimate (\ref{est:int-beta}) follows from (\ref{est:main-b}) and (\ref{est:main-c}).  

Finally, in order to verify (\ref{est:int-gamma}), we consider
the expression for $\gamma_{h,k}$, and we apply
\begin{itemize}
	\item  inequality (\ref{est:main-a}) to the term $a_{k}-a_{h}$,

	\item  inequality (\ref{est:Gamma}) to the term 
	$\Gamma_{1,k}-\Gamma_{1,h}$,

	\item  Lemma~\ref{lemma:osc-prod}, (\ref{est:main-b}), (\ref{est:main-c}) to the term $\rho_{h}^{2}(c_{h,k}-b_{h})$,

	\item  Lemma~\ref{lemma:osc-series} and (\ref{est:main-c}) to the last term (the series). 
\end{itemize}

\subsection*{Proof of Theorem 2.1} 

\paragraph{\textmd{\textit{Key estimate for quotients}}}

We prove that, if there exists $h_{0}\in J$ such that
\begin{equation}
	\limsup_{t\to +\infty}\rho_{h_{0}}(t)>0,
	\label{hp:rh>0-full}
\end{equation}
then
\begin{equation}
	\lim_{t\to +\infty}Q_{h_{0},k}(t)=1
	\quad\quad
	\forall k\in J.
	\label{lim:Q->1-full}
\end{equation}

To begin with, we observe that $\rho_{h_{0}}(t)$ is Lipschitz continuous in $[0,+\infty)$ because of (\ref{est:rk'}), and hence (\ref{hp:rh>0-full}) implies that
\begin{equation}
	\int_{0}^{+\infty}\rho_{h_{0}}^{2}(t)\,dt=+\infty.
	\label{int-div-full}
\end{equation}

Let us consider now the quotients $Q_{h_{0},k}(t)$ with $k\in J$.  We
claim that in this case equation (\ref{eqn:Q'-true}) fits in the
framework of Proposition~\ref{prop:Q1} with
$$z(t):=Q_{h_{0},k}(t),
\quad\quad
\alpha(t):=\alpha_{h_{0}}(t),
\quad\quad
\beta(t):=\beta_{h_{0},k}(t),
\quad\quad
\gamma(t):=\gamma_{h_{0},k}(t).$$

Indeed assumption (\ref{hp:Q1-a}) is exactly (\ref{int-div-full}), 
assumptions (\ref{hp:Q1-a'}) follows from (\ref{est:rk'-rk}), and 
the boundedness and semi-integrability of $\beta(t)$ and $\gamma(t)$ follow from (\ref{est:sup}), (\ref{est:int-beta}), and (\ref{est:int-gamma}).
Thus from Proposition~\ref{prop:Q1} we obtain (\ref{lim:Q->1-full}).

\paragraph{\textmd{\textit{Case where $J$ is infinite}}}

In this case we show that all components tend to 0, which establishes statement~(3).

Let us assume that this is not the case. Then there exists $h_{0}\in J$ for which (\ref{hp:rh>0-full}) holds true, and hence also  (\ref{lim:Q->1-full}) holds true. At this point, arguing exactly as in the corresponding  point in the proof of Theorem~\ref{thm:brutal}, from (\ref{hp:rh>0-full}) and (\ref{lim:Q->1-full}) we deduce that the total energy is unbounded, thus contradicting the estimate from above in (\ref{est:basic-R}).

\paragraph{\textmd{\textit{Case where $J$ is finite}}}

In this case we prove that (\ref{th:lim-rk-finite}) is true. To begin with, we observe that there exists $h_{0}\in J$ for which (\ref{hp:rh>0-full}) holds true, because otherwise the total energy would tend to~0, thus contradicting the estimate from below in (\ref{est:basic-R}). As a consequence, also (\ref{lim:Q->1-full}) holds true, and in particular the limit of $\rho_{k}(t)$ is the same for every $k\in J$, provided that this limit exists. At this point, (\ref{th:lim-rk-finite}) in equivalent to showing that
\begin{equation}
	\lim_{t\to +\infty}R(t)=\frac{4j}{2j+1},
	\label{lim-R-finite}
\end{equation}
where $j$ denotes the number of elements of $J$.

To this end we consider the equalities
$$R(t)=\rho_{h_{0}}^{2}(t)\sum_{k\in J}Q_{h_{0},k}^{2}(t),
\hspace{4em}
\sum_{k\in J}\rk^{4}(t)=
\rho_{h_{0}}^{4}(t)\sum_{k\in J}Q_{h_{0},k}^{4}(t).$$

From them we deduce that
\begin{equation}
	\sum_{k\in J}\rk^{4}(t)=
	R^{2}(t)\cdot\left(\frac{1}{j}+q(t)\right),
	\label{eqn:4-2}
\end{equation}
where
\begin{equation}
	q(t):=\left(\sum_{k\in J}Q_{h_{0},k}^{4}(t)\right)\cdot
	\left(\sum_{k\in J}Q_{h_{0},k}^{2}(t)\right)^{-2}
	-\frac{1}{j},
	\label{defn:q}
\end{equation}
hence by (\ref{lim:Q->1-full})
\begin{equation}
	\lim_{t\to +\infty}q(t)=0.
	\label{lim-q}
\end{equation} 

Going back to (\ref{eqn:R'-full}), we find that $R(t)$ solves a differential equation of the form
$$R'(t)=R(t)-\frac{2j+1}{4j}R^{2}(t)-\frac{1}{4}q(t)R^{2}(t)+\mu_1(t)+\mu_2(t),$$
where  $\mu_1(t)$ and $\mu_2(t)$ are given by (\ref{defn:f-g}). This differential equation fits in the framework of Proposition~\ref{prop:R} with
$$z(t):=R(t),
\quad\quad
z_{\infty}:=\frac{4j}{2j+1},
\quad\quad
\psi_{1}(t):=\mu_1(t)+\mu_2(t),
\quad\quad
\psi_{2}(t)=|q(t)|\cdot R^{2}(t).$$

Indeed, assumption (\ref{hp:R-psi2}) follows from (\ref{lim-q}), while
assumption (\ref{hp:R-zpos}) follows from the estimate from below in (\ref{est:basic-R}). It remains to prove that $\psi_{1}(t)$ is semi-integrable in $[0,+\infty)$. The semi-integrability of $\mu_1(t)$ is a consequence of (\ref{est:osc-R}), and the semi-integrability of $\mu_2(t)$ follows from a finite number of applications of Lemma~\ref{lemma:osc-prod} with $f(t):=\rho_{k}^{2}(t)\rho_{i}^{2}(t)$ and $g(t):=c_{i,k}(t)$ (here it is essential that the set $J$ is finite). The required assumptions of $f(t)$ and $g(t)$ follow from (\ref{est:rk}), (\ref{est:rk'}), and (\ref{est:main-c}).

At this point, Proposition~\ref{prop:R} implies (\ref{lim-R-finite}). 

\paragraph{\textmd{\textit{Asymptotic behavior of the phase}}}

It remains to prove (\ref{th:lim-tk}).  Actually we need this fact
just in the case where $J$ is finite, but the statement is true and
the proof is the same even in the general case.

Let us consider equation (\ref{eqn:tk}). From (\ref{est:Gamma}) we 
know that $\Gamma_{2,k}$ is integrable in $[0,+\infty)$. Therefore, 
(\ref{th:lim-tk}) is equivalent to showing that the function
\begin{equation}
	\left(\sum_{i=0}^{\infty}
	\rho_{i}^{2}(\tau)\sin^{2}\theta_{i}(\tau)-1\right)
	\sin\tk(\tau)\cos\tk(\tau)\
	\nonumber
\end{equation}
is semi-integrable in $[0,+\infty)$ for every $k\in J$. First of all, we write the function as
$$\sum_{i\neq k}\rho_{i}^{2}\sin^{2}\theta_{i}\sin\tk\cos\tk+
\rk^{2}\sin^{3}\tk\cos\tk-\sin\tk\cos\tk.$$

All these oscillating functions can be treated as we did many times 
before, starting from the trigonometric identities
$$\sin\tk\cos\tk=\frac{1}{2}\sin(2\tk),
\hspace{4em}
\sin^{3}\tk\cos\tk=\frac{1}{4}\sin(2\tk)-
\frac{1}{8}\sin(4\tk),$$
and 
$$\sin^{2}\theta_{i}\sin\tk\cos\tk=
\frac{1}{4}\sin(2\tk)+\frac{1}{8}\sin(2\theta_{i}-2\tk)-
\frac{1}{8}\sin(2\theta_{i}+2\tk).$$

Due to the relation $\sin x=\cos(x-\pi/2)$, we can conclude by exploiting 
the results of section~\ref{sec:oscillating} as we did in the proof 
of (\ref{est:main-a}) through (\ref{est:main-c}), and in the estimates 
of the coefficients of (\ref{eqn:Q'-true}).

\subsection*{Proof of Theorem \ref{thm:energ.}}

Let us consider again the differential equation (\ref{eqn:R'-full}) solved by $R(t)$. We prove that the uniform gap assumption (\ref{hp:diff-lambda}) implies the semi-integrability of $\mu_2(t)$ and a uniform bound on the quotients that allows to show that the series of fourth powers is negligible in the limit. At this point we can conclude by applying Proposition~\ref{prop:R}. 

\paragraph{\textmd{\textit{Estimate on $\mu_2(t)$}}}

%

We show that
\begin{equation}
	\left|\int_{t}^{s}\mu_{2}(\tau)\,d\tau\right|\leq
	M_{22}e^{-t}
	\quad\quad
	\forall s\geq t\geq 0.
	\label{est:osc-g}
\end{equation} 

Since $\mu_2(t)$ involves a double series, this requires a double application of Lemma~\ref{lemma:osc-series}. First of all, we exploit the uniform gap assumption (\ref{hp:diff-lambda}), and from (\ref{est:main-c}) we deduce that
\begin{equation}
	\left|\int_{t}^{s}c_{h,k}(\tau)\,d\tau\right|\leq M_{23}e^{-t}
	\quad\quad
	\forall s\geq t\geq 0,\ \forall h\neq k.
	\label{est:unif-cg}
\end{equation}  

Now we set
$$\delta_{k}(t):=\sum_{i\neq k}c_{i,k}(t)\rho_{i}^{2}(t),$$
and we apply Lemma~\ref{lemma:osc-series} with $f_{i}(t):=\rho_{i}^{2}(t)$ and $g_{i}(t):=c_{i,k}(t)$. The assumptions are satisfied due to (\ref{est:series-m}), (\ref{est:abc-bd}) and (\ref{est:unif-cg}). We obtain that
\begin{equation}
	\left|\int_{t}^{s}\delta_{k}(\tau)\,d\tau\right|\leq
	M_{23}e^{-t}
	\quad\quad
	\forall s\geq t\geq 0,\ \forall k\in\n.
	\label{est:delta-int}
\end{equation}
	
Moreover, from (\ref{est:abc-bd}) and (\ref{est:basic-R}) we 
obtain also that
\begin{equation}
	|\delta_{k}(t)|\leq R(t)\leq M_{4}
	\quad\quad
	\forall t\geq 0,\ \forall k\in\n.
	\label{est:delta-bd}
\end{equation}
	
Due to (\ref{est:delta-int}) and (\ref{est:delta-bd}), we can  
apply again Lemma~\ref{lemma:osc-series} with $f_{k}(t):=\rk^{2}(t)$ 
and $g_{k}(t):=\delta_{k}(t)$, and this completes the proof of 
(\ref{est:osc-g}). 
	
\paragraph{\textmd{\textit{Estimate on quotients}}}

	We claim that there exist $t_{0}\geq 0$ and $h_{0}\in J$ such that
\begin{equation}
	Q_{h_{0},k}(t)\leq 2
	\quad\quad
	\forall t\geq t_{0},\ \forall k\in J.
	\label{th:Q<2}
\end{equation}

This estimate is trivial when $k=h_{0}$, independently on $t_{0}$. 
Otherwise, we exploit equation (\ref{eqn:Q'-true}), which fits in the framework of
Proposition~\ref{prop:Q2} with
$$z(t):=Q_{h,k}(t),
\quad\quad
\alpha(t):=\alpha_{h}(t),
\quad\quad
\beta(t):=\beta_{h,k}(t),
\quad\quad
\gamma(t):=\gamma_{h,k}(t).$$  

Let us check the assumptions. Estimate (\ref{hp:Q2-bounded}) follows from (\ref{est:sup}). Estimates (\ref{hp:Q2-bc}) follow from (\ref{est:int-beta}) and (\ref{est:int-gamma}), and the constant $L_{2}$ is independent of $h$ and $k$ due to the uniform gap assumption (\ref{hp:diff-lambda}).
%
%
%
%
%
%
As a consequence, any $t_{0}\geq 0$ satisfying
(\ref{hp:Q2-t0}) is independent of $h$ and $k$, and ensures that the
implication
\begin{equation}
	Q_{h,k}(t_{0})\leq 1
	\quad\Longrightarrow\quad
	\sup_{t\geq t_{0}}Q_{h,k}(t)\leq 2
	\label{implic-Q}
\end{equation}
holds true for every $h$ and $k$ in $J$.  At this point we choose any
such $t_{0}$, and we fix the index (or one of the indices) $h_{0}\in
J$ such that
\begin{equation}
	\rho_{h_{0}}(t_{0})\geq\rk(t_{0})
	\quad\quad
	\forall k\in J.
	\nonumber
\end{equation}

Such an index exists, even when $J$ is infinite, because for every
$t\geq 0$ it turns out that $\rk(t)\to 0$ as $k\to +\infty$ due to the
square-integrability of the sequence $\rk(t)$.

This choice of $h_{0}$ implies that $Q_{h_{0},k}(t_{0})\leq 1$ for
every $k\in J$, and therefore at this point (\ref{th:Q<2}) follows
from (\ref{implic-Q}) with $h:=h_{0}$.\bigskip

\paragraph{\textmd{\textit{Conclusion}}}

To complete the proof we now observe that
$$\sum_{k\in J}\rk^{4}(t)=
\sum_{k\in J}Q_{h_{0},k}^{2}(t)\rho_{h_{0}}^{2}(t)\cdot\rk^{2}(t)\leq
4\rho_{h_{0}}^{2}(t)\cdot\sum_{k\in J}\rk^{2}(t)$$
for every $t\geq t_{0}$.
Plugging this estimate into (\ref{eqn:R'-full}) we deduce that
$$\left|R'(t)-R(t)+\frac{1}{2}R^{2}(t)-\mu_1(t)-\mu_2(t)\right|\leq
\rho_{h_{0}}^{2}(t)\cdot R(t)
\quad\quad
\forall t\geq t_{0}.$$

We are now (up to a time-translation by $t_{0}$) in the framework of
Proposition~\ref{prop:R} with
$$z(t):=R(t),
\quad\quad
z_{\infty}:=2,
\quad\quad
\psi_{1}(t):=\mu_1(t)+\mu_2(t),
\quad\quad
\psi_{2}(t):=\rho_{h_{0}}^{2}(t)\cdot R(t).$$

Indeed, the semi-integrability of $\psi_{1} $ follows from (\ref{est:osc-R}) and (\ref{est:osc-g}),
assumption (\ref{hp:R-psi2}) follows from the boundedness of $R(t)$
and the fact that $\rho_{h_{0}}(t)\to 0$ as $t\to +\infty$, and
assumption (\ref{hp:R-zpos}) follows from the estimate from below in
(\ref{est:basic-R}).

At this point, (\ref{th:lim-R-2}) is exactly the conclusion of
Proposition~\ref{prop:R}.
\qed


\subsubsection*{\centering Acknowledgments}

The first two authors are members of the Gruppo Nazionale per l'Analisi
Matematica, la Probabilit\`{a} e le loro Applicazioni (GNAMPA) of the
Istituto Nazionale di Alta Matematica (INdAM). This project was partially supported by the PRA ``Problemi di evoluzione: studio qualitativo e comportamento asintotico'' of the University of Pisa.

\label{NumeroPagine}


\begin{thebibliography}{99}

\bibitem{Bourg} \textsc{J.~Bourgain}; On the growth in time of higher Sobolev norms of smooth solutions of Hamiltonian. \emph{Internat.\ Math.\ Res.\ Notices} \textbf{6} (1996), 277--304. 

\bibitem {H-ondes} \textsc{A.~Haraux}; Comportement \`{a} l'infini pour une \'equation d'ondes non lin\'eaire dissipative (French. English summary), \emph{ C. R. Acad. Sci. Paris Sér. A-B}  \textbf{287} (1978), no. 7, A507--A509. 
\bibitem {SystD} \textsc{A.~Haraux}; \emph{Syst\`emes dynamiques dissipatifs et applications}.
Recherches en Math\'eŽmatiques Appliqu\'ees \textbf{17}, Masson, Paris, 1991. xii+132 pp. ISBN: 2-225-82283-3 .

\bibitem{Below} \textsc{A.~Haraux}; $L^p$ estimates of solutions to some non-linear wave equations in one space dimension, \emph{Int.\ J.\ Math.\ Modelling and Numerical Optimization} \textbf{1} (2009), Nos 1-2, 146--154.

\bibitem{Dieud} \textsc{A.~Haraux}; Semi-linear hyperbolic problems in bounded domains,  \emph{Math.\ Rep.}\ \textbf{3} (1987), no.~1, i--xxiv and 1--281. 

\bibitem{Open} \textsc{A.~Haraux}; Some simple problems for the next generations. Preprint arXiv:1512.06540 [math.AP] (2015).

\bibitem {H-J} \textsc{A.~Haraux, M.A. Jendoubi}; \emph{The convergence problem for dissipative autonomous systems, classical methods and recent advances}.  BCAM SpringerBriefs. Springer, Cham; BCAM Basque Center for Applied Mathematics, Bilbao, 2015. xii+142 pp., ISBN: 978-3-319-23406-9.

\bibitem{HZ} \textsc{A.~Haraux, E.~Zuazua}; Decay estimates for some semilinear damped hyperbolic problems, \emph{Arch.\ Rational Mech.\ Anal.}\ {\bf 100} (1988), no.~2, 191--206. 
 

 
\bibitem{Nakao} \textsc{M.~Nakao}; Asymptotic stability of the bounded or almost periodic solution of the wave equation with nonlinear dissipative term, \emph{J.\ Math.\ Anal.\ Appl.}\ {\bf 58} (1977), no.~2, 336--343.

\bibitem{Van-Mart} \textsc{J.~Vancostenoble, P.~Martinez};  Optimality of energy estimates for the wave equation with nonlinear boundary velocity feedbacks, \emph{SIAM J.\ Control Optim.}\ \textbf{39} (2000), no.~3, 776--797 (electronic). 
%

\end{thebibliography}
\end{document}